\documentclass[11pt,reqno,letter]{amsart}

\usepackage{amsmath}
\usepackage{amssymb}
\usepackage[dvipdfmx]{graphicx}
\usepackage{hyperref}
\usepackage[usenames]{color}
\usepackage{caption}
\usepackage{subcaption}
\usepackage{pifont}
\usepackage{bm}
\usepackage{url}
\usepackage{setspace}
\usepackage[margin=1in]{geometry}
\usepackage{subcaption}
\usepackage{xfrac}
\usepackage{appendix}

\usepackage[linesnumbered,ruled,vlined]{algorithm2e}

\usepackage{algorithmicx}

\usepackage[normalem]{ulem}

\usepackage{float}

\theoremstyle{plain}
\newtheorem{theorem}{Theorem}[section]

\newtheorem{proposition}[theorem]{Proposition}

\theoremstyle{definition}

\newtheorem{example}[theorem]{Example}
\newtheorem{assumption}[theorem]{Assumption}

\numberwithin{equation}{section}
\numberwithin{theorem}{section}
\numberwithin{table}{section}
\numberwithin{figure}{section}

\DeclareMathOperator*{\argmax}{arg\,max}
\DeclareMathOperator*{\argmin}{arg\,min}

\def\@fnsymbol#1{\ensuremath{\ifcase#1\or *\or \dagger\or \ddagger\or
   \mathsection\or \mathparagraph\or \|\or **\or \dagger\dagger
   \or \ddagger\ddagger \else\@ctrerr\fi}}

\title{A data-driven approach to beating SAA out-of-sample}

\author[Gotoh]{Jun-ya Gotoh\textsuperscript{$\dagger$}}

\author[Kim]{Michael Jong Kim\textsuperscript{$\ddagger$}}

\author[Lim]{Andrew E.B. Lim\textsuperscript{$*$}}

\dedicatory{\textsuperscript{$\dagger$}Department of Data Science for Business Innovation, Chuo University, Tokyo, Japan. Email: jgoto@kc.chuo-u.ac.jp \\ \textsuperscript{$\ddagger$}Sauder School of Business, University of British Columbia, Vancouver, Canada. Email: mike.kim@sauder.ubc.ca \\ \textsuperscript{$*$}Department of Analytics and Operations, Department of Finance, and Institute for Operations Research and Analytics, National University of Singapore, Singapore. Email: andrewlim@nus.edu.sg
}
\date{\today}

\begin{document}

\begin{abstract}
While solutions of Distributionally Robust Optimization (DRO) problems can sometimes have a higher out-of-sample expected reward than the Sample Average Approximation (SAA), there is no guarantee. In this paper, we introduce a class of Distributionally Optimistic Optimization (DOO) models, and show that it is always possible to ``beat" SAA out-of-sample if we consider not just worst-case (DRO) models but also best-case (DOO) ones.  We also show, however, that this comes at a cost:  Optimistic solutions are  more sensitive to model error than either worst-case or SAA optimizers, and hence are less robust  and calibrating the worst- or best-case model to outperform SAA may be difficult when data is limited.
\end{abstract}

\maketitle

\noindent {\bf Keywords:} Distributionally Optimistic Optimization (DOO), Distributionally Robust Optimization (DRO), Sample Average Approximation (SAA), data-driven optimization, model uncertainty, worst-case sensitivity, out-of-sample performance.


\section{Introduction}

It is well known that solutions of optimization problems calibrated from data can perform poorly out-of-sample. This occurs due to errors in both the modeling assumptions and the estimation of model parameters and has motivated various optimization models which account for misspecification in the in-sample model. Distributionally Robust Optimization (DRO), where the decision maker optimizes against worst-case perturbations from the nominal (in-sample) distribution, is one such approach. Similar models have been introduced in a number of communities \cite{BN2,DGKF,EL,Hansen3,PJD}.


The performance of DRO is often evaluated by comparing its out-of-sample expected reward to that of the SAA optimizer. This is often done experimentally, and recent papers characterize finite-sample and asymptotic properties of DRO solutions and the associated expected reward \cite{bgk-SAA,DN,DGN16,GKL-cal,KENS}.

It has  been observed that solutions of DRO and other worst-case models can sometimes have a higher out-of-sample expected reward than the Sample Average Approximation (SAA) optimizer \cite{AP,BLSZZ2013,GKL-cal,KL}.
However, there is no guarantee, and it is of interest whether a decision that beats SAA out-of-sample can be found when DRO is unable to.
In this paper, we introduce a class of {\it Distributionally Optimistic Optimization} (DOO) problems, where nature works together with the decision maker to optimize the in-sample reward, and show under relatively mild assumptions that it is always possible to ``beat" the SAA optimizer out-of-sample if we consider best-case (DOO) and worst-case (DRO) decisions.

As tempting as this may sound,  there is a catch. While the out-of-sample expected reward might be larger, the expected reward under a best-case optimizer is always more sensitive to worst-case perturbations from the nominal model than the SAA optimizer, and hence is less robust. More generally, worst and best-case optimizers make a tradeoff between (in-sample) expected reward and  sensitivity to model error. We potentially miss the benefits of sensitivity reduction if we only focus on the mean reward, and could even make things worse if this leads  to choosing a  best-case solution.

In contrast to worst-case problems, the literature on DOO is small, though recent activity suggests an uptick in interest that is generally centered around the concern that worst-case optimization is ``too pessimistic." One disadvantage of optimistic optimization, however, is that convexity can be lost, making it potentially more difficult to solve (computationally) than its worst-case cousin. 
Non-convexity is considered in \cite{NTM} where optimistic optimization problems are shown to be related to non-convex regularizers used in machine learning. 
A majorization minimization algorithm is employed in \cite{NTM} to efficiently obtain a critical point of their nonconvex optimization problem by making use of its ``difference of convex" structure. 
Also related is \cite{NSYKW}, which uses optimistic optimization over probability distributions to construct a nonparametric approximation of the likelihood function that can be used in Bayesian statistics. 
The  paper \cite{SZ} uses DOO to solve the Trust Region Policy Optimization problem in the context of Reinforcement Learning, while optimistic optimization is also proposed in \cite{CR} for online contextual decision making. Finally, a recent paper \cite{CR} considers best-case optimization in stochastic optimization with model uncertainty using the ``Rockafellian" framework for perturbation analysis. Best-case problems are also used in the empirical likelihood approach to generating confidence intervals of the optimal expected reward under the population distribution \cite{DGN16,Lam1,LamZhou,WGY}.  We also note that the National Science Foundation recently awarded a grant \cite{Xie} on ``Favorable optimization under distributional distortions". These papers do not consider the out-of-sample performance of solutions of best-case problems, nor the impact on the sensitivity of the out-of-sample expected reward to errors in the in-sample model.

A recent paper \cite{Lam2} shows that it is impossible, {\it asymptotically}, for a large class of data-driven problems, including DRO and DOO and popular regularization techniques, to beat the out-of-sample expected reward of SAA. While this appears to directly contradict our results, there is actually no inconsistency. The main difference is that \cite{Lam2} considers the large data limit, whereas our results apply to data sets of moderate size. Taken together, it is always possible to beat SAA using DRO/DOO with a finite data set (us), but impossible in the limit \cite{Lam2}.
We discuss \cite{Lam2} in more detail at the end of this paper.

\subsection*{Organization}
We introduce the {\it Distributionally Optimistic Optimization} (DOO) and {\it Distributionally Robust Optimization} (DRO) models in Section \ref{sect:setup}. For concreteness, we adopt a penalty framework with smooth $\phi$-divergence as the penalty function.
We show in Section \ref{sect:solutions} that the family of worst- and best-case solutions is continuously differentiable in a neighborhood of the SAA solution and characterize their distributional properties, that it is generally possible to find a DRO or DOO solution with a higher out-of-sample expected reward than SAA in Section \ref{sect:out-of-sample}, but that a best-case solution is  always less robust than SAA in Section \ref{sect:wcs}. The key ideas are illustrated using a data-driven inventory problem.

\section{Setup} \label{sect:setup}

Let $f:{\mathbb R}^d\times{\mathbb R}^l \rightarrow {\mathbb R}$, and $Y$ be an  $\mathbb{R}^{l}$-valued random vector with population distribution $\mathbb{P}$.
Consider the problem
\begin{eqnarray}
\max_x {\mathbb E}_{\mathbb P} \big[f(x,\,Y)\big]. \label{eq:pop}
\end{eqnarray}
Let
\begin{equation}
x^{\star}(0) := \arg \max_x \,\Big\{ \mathbb{E}_{\mathbb P}\left[f(x,\,Y)\right] \Big\}
\label{eq:opt_dgm}
\end{equation}
denote the solution of \eqref{eq:pop}. We assume the following.
\begin{assumption} \label{ass:f}
The function $f(x,\,Y)$ and random vector $Y\sim {\mathbb P}$ are such that
\begin{itemize}
\item $f(x,\,Y)$ is strictly concave and twice continuously differentiable in $x\in\mathbb{R}^d$ for ${\mathbb P}$-almost surely every $Y\in\mathbb{R}^{l}$;
\item for each fixed $x\in\mathbb{R}^{d}$, the mappings $y\mapsto f(x,y), \nabla_{x}f(x,y), \nabla^{2}_{x}f(x,y)$, $y\in\mathbb{R}^{l}$, are measurable and all moments of the random variables $f(x,Y)$, $\nabla_{x}f(x,Y)$, $\nabla^{2}_{x}f(x,Y)$ exist;
\item there exists a solution $x^\star(0)$ of  \eqref{eq:pop}.
\end{itemize}
\end{assumption}

In many applications, we do not know the population distribution ${\mathbb P}$ but only have independent and identically distributed (iid) samples $Y_1,\cdots,\,Y_n$ drawn from $\mathbb P$. This naturally leads to replacing \eqref{eq:pop} with a Sample Average Approximation (SAA):
\begin{eqnarray}
\max_x \; {\mathbb E}_{{\mathbb P}_n} \left[f(x, \, Y)\right] \equiv \sum_{i=1}^n p_i ^n f(x,\,Y_i). \label{eq:SAA}
\end{eqnarray}
Here ${\mathbb P}_n = [p_1^{n},\cdots,\,p_n^{n}]$ is the empirical distribution; we assume without loss of generality that $p^n_i>0$ for all $i$. We denote the solution of the in-sample problem \eqref{eq:SAA} by
\begin{eqnarray*}
x_n(0) := \argmax_x \sum_{i=1}^n p_i^n f(x,\,Y_i).
\end{eqnarray*}
It is well known that the SAA solution $x_n(0)$ may not perform well out-of-sample. This has motivated worst-case versions of SAA, called Distributionally Robust Optimization (DRO), where the decision is chosen to maximize the expected reward under worst-case perturbations of the probability distribution (${\mathbb Q}$) from the empirical distribution ${\mathbb P}_n$:
\begin{eqnarray}
\max_x \min_{\mathbb Q}\sum_{i=1}^n q_i f(x,\,Y_i) + \frac{1}{\delta}{\mathcal H}_\phi({\mathbb Q}\,|\,{\mathbb P}_n),\; \delta>0 \label{eq:DRO}
\end{eqnarray}
where
\begin{eqnarray}
\mathcal{H}_\phi(\mathbb{Q} \,|\, {\mathbb{P}}_n):=
\sum\limits_{i=1}^n {p}^n_i\phi\left(\frac{q_i}{{p}^{n}_i}\right),&\sum\limits_{i=1}^n q_i=1, q_i\geq 0
\label{phi-div}
\end{eqnarray}
is the $\phi$-divergence of ${\mathbb Q}=[q_1,\cdots,\,q_n]$ relative to ${\mathbb P}_n = [p_1^{n},\cdots,\,p_n^{n}]$.
We assume throughout that $\phi:{\mathbb R}\rightarrow{\mathbb R}\cup\{+\infty\}$ is a convex lower semi-continuous function such that $\phi(z)\geq\phi(1)=0$ for $z\geq 0$ and $\phi(z)=+\infty$ for $z<0$.

\subsection*{Distributionally optimistic optmization (DOO)}
We now consider an optimistic version of \eqref{eq:SAA} where nature {\it works together} with the decision maker to optimize the expected reward. The Distributionally Optimistic Optimization (DOO) model is
\begin{eqnarray}
\max_x \max_{\mathbb Q}\sum_{i=1}^n q_i f(x,\,Y_i) + \frac{1}{\delta}{\mathcal H}_\phi({\mathbb Q}\,|\,{\mathbb P}_n), \; \delta<0.  \label{eq:OPT}
\end{eqnarray}
Here, nature selects $\mathbb Q$ to maximize the expected reward and incurs a (negative) penalty for deviating from the nominal ${\mathbb P}_n$; the parameter $\delta$ controls the size of the deviations.
The decision maker accepts that ${\mathbb P}_n$ may be misspecified, and makes a decision that maximizes the expected reward if nature is cooperative.
DOO may be of interest to applications where the decision maker is looking for upside opportunities and does not wish to be conservative in the face of model uncertainty.
We also note that one of the likelihood approximations in \cite{NSYKW} has the form of the inner problem in \eqref{eq:OPT}. They do not consider optimization, though an optimistic version of Maximum Likelihood Estimation is natural step in this direction and an example of \eqref{eq:OPT}.
 Though our results provide insights about general properties of DOO and its solutions, applications of DOO are not the focus of this paper.


It is convenient to introduce population versions of the worst-case and best-case problems. If
\begin{eqnarray}
\mathcal{H}_\phi(\mathbb{Q} \,|\, {\mathbb{P}}):=
{\mathbb E}_{\mathbb P}\Big[\phi\Big(\frac{\mathrm{d}{\mathbb Q}}{\mathrm{d}{\mathbb P}} \Big)\Big]
\label{phi-div_pop}
\end{eqnarray}
denotes $\phi$-divergence,  where $\frac{\mathrm{d}\mathbb Q}{\mathrm{d}\mathbb P}$ is the Radon-Nikodym derivative (likelihood ratio) of $\mathbb Q$ with respect to $\mathbb P$, then the population version of the worst-case problem is
\begin{eqnarray}
\max_x \min_{\mathbb Q} {\mathbb E}_{\mathbb Q} [f(x,\,Y)] + \frac{1}{\delta}{\mathcal H}_\phi({\mathbb Q}\,|\,{\mathbb P}),\; \delta>0 \label{eq:DRO_pop}
\end{eqnarray}
and
\begin{eqnarray}
\max_x \max_{\mathbb Q}{\mathbb E}_{\mathbb Q} [f(x,\,Y)] + \frac{1}{\delta}{\mathcal H}_\phi({\mathbb Q}\,|\,{\mathbb P}), \; \delta<0  \label{eq:OPT_pop}
\end{eqnarray}
is the best-case problem.

We denote the solutions of \eqref{eq:opt_dgm}, \eqref{eq:DRO_pop} and \eqref{eq:OPT_pop} by $x^\star(\delta)$ with $\delta=0$, $\delta>0$ and $\delta<0$, respectively, and $x_n(\delta)$ for the sample versions \eqref{eq:SAA}, \eqref{eq:DRO} and \eqref{eq:OPT}. It will be shown in Section \ref{sect:solutions} that $x^\star(\delta)$ and $x_n(\delta)$ are continuously differentiable in a neighborhood of $\delta=0$.  We consider the case where the decision variable $x$ is unconstrained as this allows us to directly connect changes to the SAA solution from DRO/DOO to out-of-sample performance. We defer analysis of the constrained case to future work.

\subsection*{Dual characterization of worst/best-case objective}

Let
\begin{eqnarray*}
\phi^*(\zeta) := \max_z\Big\{\zeta z - \phi(z)\Big\}
\end{eqnarray*}
denote the convex conjugate of $\phi(z)$.
We have the following dual representation of the worst-case problems \eqref{eq:DRO} and \eqref{eq:DRO_pop}, which is well known and can be established using convex duality. We state it without proof though the interested reader can see \cite{GKL}.
\begin{proposition}[Dual characterization for DRO]\label{prop:DRO_dualrep}
Suppose that $\phi:{\mathbb R}\rightarrow{\mathbb R}\cup\{+\infty\}$ is a convex lower-semicontinuous function such that $\phi(z)\geq\phi(1)=0$ for $z\geq 0$ and $\phi(z)=+\infty$ for $z<0$. If $\delta>0$, then
\begin{eqnarray}
\min_{\mathbb Q} \Big\{ {\mathbb E}_{\mathbb Q}[f(x,\,Y)] + \frac{1}{\delta}{\mathcal H}_\phi({\mathbb Q}\,|\,{\mathbb P})\Big\}= \max_c  \Big\{-\frac{1}{\delta}{\mathbb E}_{\mathbb P}\Big[\phi^*\Big(-\delta\big[f(x,\,Y)+ c\big]\Big)\Big]-c\Big\} \label{eq:DRO_dual_pop}
\end{eqnarray}
and
\begin{eqnarray}
\min_{\mathbb Q}\Big\{\sum_{i=1}^n q_i f(x,\,Y_i) + \frac{1}{\delta}{\mathcal H}_\phi({\mathbb Q}\,|\,{\mathbb P}_n)\Big\} = \max_c  \Big\{-\frac{1}{\delta}\sum_{i=1}^np_i^n \phi^*\Big(-\delta\big[f(x,\,Y_i)+ c\big]\Big)-c\Big\}. \label{eq:DRO_dual}
\end{eqnarray}
\end{proposition}
Similarly,  we can derive the dual representation for the optimistic problems  \eqref{eq:OPT} and \eqref{eq:OPT_pop}.
\begin{proposition}[Dual characterization for DOO]\label{prop:OPT_dualrep}
Suppose that $\phi:{\mathbb R}\rightarrow{\mathbb R}\cup\{+\infty\}$ is a convex lower-semicontinuous function such that $\phi(z)\geq\phi(1)=0$ for $z\geq 0$ and $\phi(z)=+\infty$ for $z<0$. If $\delta<0$, then
\begin{eqnarray}
\max_{\mathbb Q} \Big\{ {\mathbb E}_{\mathbb Q}[f(x,\,Y)] + \frac{1}{\delta}{\mathcal H}_\phi({\mathbb Q}\,|\,{\mathbb P})\Big\} = \min_c  \Big\{-\frac{1}{\delta}{\mathbb E}_{\mathbb P}\Big[\phi^*\Big(-\delta\big[f(x,\,Y)+ c\big]\Big)\Big]-c\Big\} \label{eq:OPT_dual_pop}
\end{eqnarray}
and
\begin{eqnarray}
\max_{\mathbb Q}\sum_{i=1}^n q_i f(x,\,Y_i) +\frac{1}{\delta}{\mathcal H}_\phi({\mathbb Q}\,|\,{\mathbb P}_n)= \min_c \Big\{-\frac{1}{\delta}\sum_{i=1}^np_i^n \phi^*\Big(-\delta\big[f(x,\,Y_i)+ c\big]\Big)-c\Big\}.\label{eq:OPT_dual}
\end{eqnarray}
\end{proposition}

\section{Characterization of the solution} \label{sect:solutions}

\subsection{In-sample problems}
It follows from Proposition \ref{prop:DRO_dualrep}  that the worst-case problems \eqref{eq:DRO} and \eqref{eq:DRO_pop} are equivalent to
\begin{align}
\max_x \max_c  \Big\{-\frac{1}{\delta}{\mathbb E}_{\mathbb P}\Big[\phi^*\Big(-\delta\big[f(x,\,Y)+ c\big]\Big)\Big]-c\Big\},\;\delta>0 \label{eq:opt_DRO}\\[5pt]
\max_x \max_c  \Big\{-\frac{1}{\delta}\sum_{i=1}^np_i^n\phi^*\Big(-\delta\big[f(x,\,Y_i)+ c\big]\Big)-c\Big\},\;\delta>0 \nonumber
\end{align}
while by Proposition \ref{prop:OPT_dualrep}, the best case problems \eqref{eq:OPT} and \eqref{eq:OPT_pop} become
\begin{align}
\max_x \min_c  \Big\{-\frac{1}{\delta}{\mathbb E}_{\mathbb P}\Big[\phi^*\Big(-\delta\big[f(x,\,Y)+ c\big]\Big)\Big]-c\Big\}, \; \delta<0 \label{eq:opt_DOO}\\[5pt]
\max_x \min_c \Big\{-\frac{1}{\delta}\sum_{i=1}^np_i^n \phi^*\Big(-\delta\big[f(x,\,Y_i)+ c\big]\Big)-c\Big\},\; \delta<0.\nonumber
\end{align}
We study properties of the solution of these problems using the first-order conditions. The assumptions about the function $\phi(z)$ need to be strengthened.
\begin{assumption}
\label{ass:phi}
$\phi:{\mathbb R}\rightarrow{\mathbb R}\cup\{+\infty\}$ is a
convex lower-semicontinuous function such that $\phi(z)\geq\phi(1)=0$ for $z\geq 0$, $\phi(z)=+\infty$ for $z<0$, and is twice continuously differentiable around $z=1$ with $\phi''(1)>0$.
\end{assumption}
The first order conditions for the sample and population problems in  \eqref{eq:opt_DRO} and \eqref{eq:opt_DOO} can be written in the form
\begin{eqnarray}
{\mathbb E}_{{\mathbb P}_n}\big[\psi(x,\,c, \, Y)\big] = \left[\begin{array}{cc}0 \\ 0\end{array}\right]
\label{eq:FOC_emp}
\end{eqnarray}
and
\begin{eqnarray}
{\mathbb E}_{{\mathbb P}}\big[\psi(x,\,c,\, Y)\big] = \left[\begin{array}{cc}0 \\ 0\end{array}\right],
\label{eq:FOC_pop}
\end{eqnarray}
respectively, where $\delta>0$ for DRO and $\delta<0$ for DOO, and
\begin{eqnarray}
\psi(x,\,c, \, Y) \equiv  \left[\begin{array}{c} \psi_1(x, \, c, \, Y) \\ \psi_2(x, \, c, \, Y)\end{array}\right] :=
\left[
\begin{array}{c}
\displaystyle [\phi^*]'\big(-\delta\big[f(x,\,Y)+c\big]\big)
\nabla_xf(x,\,Y)  \\[5pt]
\displaystyle
-\frac{\phi''(1)}{\delta} \Big\{[\phi^*]'\big(-\delta\left(f(x,\,Y)+c\right)\big) -1\Big\}
\end{array}
\right].
\label{eq:psi}
\end{eqnarray}
(see Appendix \ref{app-psi} for a derivation). Since 
\begin{eqnarray*}
\phi^*(\zeta) = \zeta + \frac{1}{2!}\frac{\zeta^2}{\phi''(1)} +o(\zeta^2)
\end{eqnarray*}
(see Theorem 3.2 in \cite{GKL}), it follows that
\begin{eqnarray*}
[\phi^*]'(\zeta) = 1 + \frac{\zeta}{\phi''(1)}  + o(\zeta),
\end{eqnarray*}
and
\begin{eqnarray*}
\psi(\,x, \, c, \, Y)=\left[\begin{array}{c} \nabla_x f(x,Y) -\frac{\delta}{\phi''(1)} \Big(f(x,Y)+c\Big)\nabla_x f(x,Y)+ o(\delta) \\[10pt] 
f(x,Y)+c + O(\delta) \end{array}\right].
\end{eqnarray*}

To ease notation, we suppress the dependence on $\delta$ in the function $\psi(x,\,c,\,Y)$; it should be clear from the context whether we are talking about the worst-case  ($\delta>0$) or the best-case  ($\delta<0$) problem.

Let $(x_n(\delta),\,c_n(\delta))$ and $(x^\star(\delta),\,c^\star(\delta))$ denote the solutions of the in-sample \eqref{eq:FOC_emp} and population \eqref{eq:FOC_pop} problems, respectively. The following result shows that the family of solutions parameterized by $\delta$ exists and is continuously differentiable in a neighborhood of $\delta=0$. There is a similar result in \cite{GKL-cal} though the focus there is limited to DRO problems ($\delta\geq 0$). Proposition \ref{prop:solution_differentiability} shows that the continuation to negative values of $\delta$ are solutions of DOO problems. The proof is in the Appendix.
\begin{proposition} \label{prop:solution_differentiability}
Suppose that $f(x,\,Y)$ satisfies Assumption \ref{ass:f} and $\phi(z)$ satisfies Assumption \ref{ass:phi}. Then there is an open neighborhood of $\delta=0$ where the solutions
$(x_n(\delta),\,c_n(\delta))$ and $(x^\star(\delta),\,c^\star(\delta))$ of \eqref{eq:FOC_emp}  and \eqref{eq:FOC_pop}, respectively, exist and are continuously differentiable.
 In particular,
\begin{align}
\left\{
\begin{array}{l} x_n(\delta)  = x_n(0) + \pi_n \delta + o(\delta),\\ [5pt]
c_n(\delta)  = - {\mathbb E}_{{\mathbb P}_n} [f(x_n(0),\,Y)]  + O(\delta),
\end{array}\right.
\label{eq:insample_rob_asymp_bias}
\end{align}
where\footnote{If $\phi(z)$ is three times continuously differentiable, it can be shown that \begin{eqnarray*}c_n(\delta) = - {\mathbb E}_{{\mathbb P}_n} [f(x_n(0),\,Y)] - \frac{\delta}{2}\frac{\phi'''(1)}{\left[\phi''(1)\right]^2} \mathbb{V}_{{\mathbb P}_n}\big[f(x_n(0),\,Y)\big] + o(\delta).\end{eqnarray*} However, the first-order term is not needed in our analysis.}
\begin{align}
\pi_n  := \frac{1}{\phi^{''}(1)} \, \Big({\mathbb{E}_{{\mathbb{P}}_n}[\nabla^2_x f(x_n(0),\,Y)]}\Big)^{-1} \mathrm{Cov}_{{\mathbb{P}_n}}\Big[{\nabla_x f(x_n(0),\,Y)}, \, f(x_n(0),\,Y)\Big],
\label{eq:insample_pi}
\end{align}
and
\begin{align}
 \left\{\begin{array}{l}x^\star(\delta)  = x^\star(0) + \pi \delta + o(\delta), \\ [5pt]
c^\star(\delta)  = - {\mathbb E}_{{\mathbb P}} [f(x^\star(0),\,Y)]  + O(\delta),
\end{array}\right.
\label{eq:pop_rob_asymp_bias}
\end{align}
where
\begin{align}
\pi := \frac{1}{\phi^{''}(1)} \, \Big({\mathbb{E}_{{\mathbb{P}}}[\nabla^2_x f(x^\star(0),\,Y)]}\Big)^{-1} \mathrm{Cov}_{{\mathbb{P}}}\Big[{\nabla_x f(x^\star(0),\,Y)}, \, f(x^\star(0),\,Y)\Big].
\label{eq:pop_pi}
\end{align}
\end{proposition}
Proposition \ref{prop:solution_differentiability} shows that worst- and best-case optimization adds a bias in the direction $\pi_n$ to the SAA maximizer $x_n(0)$ (and similarly for $x^\star(\delta)$).


Clearly, the in-sample expected reward under the robust optimizer is smaller than that of the empirical optimizer, no matter the sign of $\delta$.
When  applied out-of-sample, however, this might not be the case because the solutions are random variables under the population distribution $\mathbb P$. To compare the out-of-sample reward of the best and worst-case solutions and the SAA optimizer, we need to understand the impact of data variability on the solution.

\subsection{Statistical properties of the SAA solution} \label{sec:statistical}

The following regularity assumptions on  $f(x, Y)$ will be needed to characterize the statistical properties of $x_n(0)$.

\begin{assumption} \label{ass:regf}
\begin{enumerate}
\item The function 
$f(x, Y)$ is three times continuously differentiable in $x\in{\mathbb R}^d$ for ${\mathbb P}$-almost all $Y\in{\mathbb R}^l$;
\item As $n\rightarrow\infty$
\begin{eqnarray*}
\sup_{x}\Big|{\mathbb E}_{{\mathbb P}_n}[\nabla_x f(x,\,Y)]-{\mathbb E}_{{\mathbb P}}[\nabla_x f(x, \,Y)]\Big| \overset{P}{\longrightarrow}  0
\end{eqnarray*}
\item for every $\epsilon>0$
\begin{eqnarray*}
\inf_{x: \left|x-x^\star(0)\right| \geq \epsilon} \Big|{\mathbb E}_{\mathbb P} [\nabla_x f(x,\,Y)]\Big| > 0 = \Big|{\mathbb E}_{\mathbb P} [\nabla_x f(x^\star(0),\,Y)]\Big|.
\end{eqnarray*}
\item The $2^{nd}$ and $3^{rd}$ order derivatives of $f (x, \, Y)$ exist for any $x$ in a neighborhood of $x^\star(0)$ and
\begin{eqnarray*}
{\mathbb E}_{\mathbb P} \big|\nabla_x^2f(x, Y)\big|^2 & < & \infty, \\
{\mathbb E}_{\mathbb P} \Big|\nabla_x^2 \Big(\frac{\partial f}{\partial x_i}(x, \, Y)\Big)\Big|^2 & < & \infty, \; i=1,\cdots, d
\end{eqnarray*}
where $\nabla_x^2 \Big(\frac{\partial f}{\partial x_i}(x, \, Y)\Big)$ is the ${\mathbb R}^{d\times d}$-valued Hessian of the real-valued function $\frac{\partial f}{\partial x_i}(x, \, Y)$.
\item For any $x$ in some neighborhood of $x^\star(0)$
\begin{eqnarray*}
\left({\mathbb E}_{{\mathbb P}_n}\big[\nabla_x^2 f(x, Y)\big]\right)^{-1} = O_P(1)
\end{eqnarray*}
\item For any $x$ in some neighborhood of $x^\star(0)$ and random variable $M$ such that ${\mathbb E}_{\mathbb P}|M| < \infty$,
\begin{eqnarray*}
{\mathbb E}_{\mathbb P} \Big|\nabla_x^2 f(x,  Y)-\nabla_x^2f(x^\star(0),  Y)\Big| & < & M |x-x^\star(0)|, \\ [5pt]
{\mathbb E}_{\mathbb P} \Big|\nabla_x^2 \Big(\frac{\partial f}{\partial x_i}(x, \, Y)\Big)- \nabla_x^2 \Big(\frac{\partial f}{\partial x_i} (x^\star(0), \, Y)\Big)\Big| & < & M |x-x^\star(0)|.
\end{eqnarray*}
\end{enumerate}
\end{assumption}


We now characterize the statistical properties of the SAA solution.
Let
\begin{eqnarray*}
{A}(0) & := &  {\mathbb E}_{\mathbb{P}}\big[\nabla^2_xf(x^{\star}(0),\, Y)\big] \in {\mathbb R}^{d\times d},\\
{B}(0)  & := &  {\mathbb E}_{\mathbb P}\big[\nabla_x f(x^\star(0),\, Y)\,\nabla_x f (x^\star(0),\, Y)'\big] \in {\mathbb R}^{d \times d},
\end{eqnarray*}
and random vectors
\begin{eqnarray*}
{W}_n(0) & := & -A(0)^{-1}\frac{1}{\sqrt{n}}\sum_{i=1}^n\nabla_x f(x^\star(0), \, Y_i), \\
{V}_n(0) & := & H_n(0) W_n(0) - I_n(0),
\end{eqnarray*}
where
\begin{eqnarray*}
{H}_n(0) & := & -\frac{1}{\sqrt{n}}{A}(0)^{-1}\sum_{i=1}^n\Big\{ \nabla^2_xf(x^{\star}(0),\, Y_i) - {\mathbb E}_{\mathbb P}\big[\nabla^2_x f(x^{\star}(0),\, Y)\big]\Big\} \\
{I}_n(0) & := & {A}(0)^{-1}
\left[\begin{array}{c}
{W}_n(0)'{\mathbb E}_{\mathbb P}\left[\nabla^2_x\big(\frac{\partial f}{\partial {x_1}} (x^\star(0),\, \, Y)\big) \right]{W}_n(0)  \\
\vdots \\
{W}_n(0)'{\mathbb E}_{\mathbb P}\left[\nabla^2_x\big(\frac{\partial f}{\partial x_d}(x^\star(0),\, Y)\big)\right]{W}_n(0)
\end{array}\right].
\end{eqnarray*}
By condition (1) of Assumption \ref{ass:regf}, $\nabla^2_x\big(\frac{\partial f}{\partial x_i}(x^\star(0),\, Y)\big)$ is guaranteed to exist. 

The following result gives the statistical properties of the solution of the SAA problem.
\begin{proposition} \label{prop:SAA}
Suppose that data $\{Y_1,\cdots,\,Y_n\}$ is drawn iid from $\mathbb P$ and  $f(x,\,Y)$ satisfies Assumptions \ref{ass:f} and \eqref{ass:regf}.
Then there is a unique solution
$x^\star(0)$ of the first-order conditions ${\mathbb E}_{{\mathbb P}_n}[f(x, Y)]=0$ for the SAA problem \eqref{eq:SAA}, the matrix ${A}(0)$ is invertible,  $x_n(0) \overset{P}{\longrightarrow} x^\star(0)$ as $n\rightarrow\infty$, and 
\begin{eqnarray}
x_n(0) = x^\star(0) + \frac{1}{\sqrt{n}} W_n(0) + \frac{1}{n} V_n(0) + o_P(n^{-3/2})
\label{eq:x0}
\end{eqnarray}
where ${W}_n(0)$ has mean $0$ and covariance matrix
\begin{eqnarray*}
\xi(0) & = & A(0)^{-1} B(0) {A(0)^{-1}}' \\
 & = &  {\mathbb V}_{\mathbb P}\Big[\Big({\mathbb E}\big[\nabla^2_x f(x^{\star}(0),\,Y)\big]\Big)^{-1}{\nabla_x f(x^{\star}(0),\,Y)}\Big],
\end{eqnarray*}
and ${V}_n(0)$ is $O_P(1)$ with mean
\begin{eqnarray}
\lefteqn{\overline{V}(0) = {\mathbb E}_{\mathbb P} \Big[{A}(0)^{-1} \nabla^2_xf(x^{\star}(0),\, Y) {A}(0)^{-1}\nabla_x f (x^{\star}(0),\, Y)\Big]} \nonumber \\
& &
- {A}(0)^{-1}
\left[\begin{array}{c}
{\rm tr} \big\{\xi(0) \, {\mathbb E}_{\mathbb P}\big[\nabla^2_x\frac{\partial f}{\partial x_1}(x^\star(0),\,Y) \big] \big\} \\
\vdots  \\
{\rm tr} \big\{\xi(0)\,{\mathbb E}_{\mathbb P}\big[\nabla^2_x\frac{\partial f}{\partial x_d} (x^\star(0),\, Y)\big]\big\}
\end{array}\right].
\label{eq:V-bar}
\end{eqnarray}
\end{proposition}
\begin{proof}
Uniqueness of $x^\star(0)$ was established in Proposition \ref{prop:solution_differentiability}, while consistency of $x_n(0)$ follows from  Theorem 5.9 in \cite{vdV}. Equation \eqref{eq:x0} is obtained by applying Lemma 2.1 in \cite{KR} or Lemma 3.1 from \cite{RSU} to the first-order conditions of the SAA problem \eqref{eq:SAA}.  Invertability of $A(0)$ follows from Assumption \ref{ass:f}.
\end{proof}

Equation \eqref{eq:x0} shows that the SAA solution has a bias $\frac{1}{n}\overline{V}(0)$ relative  to the solution $x^\star(0)$ of the population problem.

\subsection{Statistical properties of the DRO/DOO solution}\label{sec:statistical-DRO}
The following  assumptions for $\psi(x, c, Y)$ are a direct parallel of Assumption \ref{ass:regf} for $f(x, Y)$. 
\begin{assumption} \label{ass:reg}
\begin{enumerate}
\item The function 
$f(x, Y)$ is three times continuously differentiable in $x\in{\mathbb R}^d$ for ${\mathbb P}$-almost surely $Y\in{\mathbb R}^l$, and $\phi(z)$ is three times continuously differentiable in $z\in{\mathbb R}$, and hence, $\phi^*(\zeta)$ is three times continuously differentiable in $\zeta$.
\end{enumerate}
Let $\psi(x, c, Y)$ be defined by \eqref{eq:psi}. The parameter $\delta$ in $\psi(x, c, Y)$ is such that
\begin{enumerate}
\item[(2)] As $n\rightarrow\infty$
\begin{eqnarray*}
\sup_{(x,\,c)}\Big|{\mathbb E}_{{\mathbb P}_n}[\psi(x,\,c,\,Y)]-{\mathbb E}_{{\mathbb P}}[\psi(x,\,c,\, \,Y)]\Big| \overset{P}{\longrightarrow}  0
\end{eqnarray*}
\item[(3)] for every $\epsilon>0$
\begin{eqnarray*}
\inf_{(x,\,c): |(x,\,c)-(x^\star(\delta),\,c^\star(\delta))|\geq \epsilon} \Big|{\mathbb E}_{\mathbb P} [\psi(x,\,c,\,Y)]\Big| > 0 = \Big|{\mathbb E}_{\mathbb P} [\psi(x^\star(\delta),\,c^\star(\delta),\,Y)]\Big|.
\end{eqnarray*}
\item[(4)] The $1^{st}$ and $2^{nd}$ order derivatives of $\psi(x, c, Y)$ exist in a neighborhood of $(x^\star(\delta), c^\star(\delta))$ and
\begin{eqnarray*}
{\mathbb E}_{\mathbb P} \big|\nabla_{(x, c)}\psi(x, c, Y)\big|^2 & < & \infty, \\
{\mathbb E}_{\mathbb P} \big|\nabla_{(x, c)}^2 \psi^{(i)}(x, c, Y)\big|^2 & < & \infty, \; i=1,\cdots, d+1
\end{eqnarray*}
where $\psi^{(i)}$ be the $i^{th}$ component of $\psi = [\psi^{(1)},\cdots,\psi^{(d)}, \psi^{(d+1)}]'$. 
\item[(5)] For some neighborhood of $(x^\star(\delta), c^\star(\delta))$
\begin{eqnarray*}
\left({\mathbb E}_{{\mathbb P}_n}\big[\nabla_{(x, c)}\psi(x, c, Y_i)\big]\right)^{-1} = O_P(1)
\end{eqnarray*}
\item[(6)] For some neighborhood of $(x^\star(0), c^\star(0))$ and random variable $M$ such that ${\mathbb E}_{\mathbb P}|M| < \infty$,
\begin{eqnarray*}
{\mathbb E}_{\mathbb P} \Big|\nabla_{(x, c)}\psi(x, c, Y)-\nabla_{(x, c)}\psi(x^\star(0), c^\star(0), Y)\Big| & < & M\Big(|x-x^\star(0)| + |c-c^\star(0)|\Big), \\ [5pt]
{\mathbb E}_{\mathbb P} \Big|\nabla_{(x, c)}^2 \psi^{(i)}(x, c, Y)- \nabla_{(x, c)}^2 \psi^{(i)}(x^\star(0), c^\star(0), Y)\Big| & < & M\Big(|x-x^\star(0)| + |c-c^\star(0)|\Big).
\end{eqnarray*}
\end{enumerate}
\end{assumption}
Condition (1) guarantees that the Hessian $\nabla_{(x, c)}^2 \psi^{(i)}(x, c, Y)$ exists. The third condition implies that $\delta$ is such that the first order conditions \eqref{eq:FOC_pop} for the population problem has a unique solution. Under Assumptions \ref{ass:f} and \ref{ass:phi}, this is satisfied for every $\delta\geq 0$
(i.e. for SAA and the worst-case problems), and  for all $\delta$ in an open neighborhood of $0$, which includes negative values.

The following result characterizes the statistical properties of the DRO/DOO solution $x_n(\delta)$.
\begin{proposition} \label{prop:distribution-of-solution}
Suppose that data $\{Y_1,\cdots,\,Y_n\}$ is drawn iid from $\mathbb P$, $f(x,\,Y)$ satisfies Assumptions \ref{ass:f} and \ref{ass:regf},  $\phi(z)$ satisfies Assumption \ref{ass:phi}, and 
 $\delta$ is such that Assumption \ref{ass:reg} holds.
Then there is a unique solution $(x^\star(\delta),\,c^\star(\delta))$ of the first-order conditions \eqref{eq:FOC_pop}, $(x_n(\delta),\,c_n(\delta)) \overset{P}{\longrightarrow} (x^\star(\delta),\,c^\star(\delta))$ as $n \rightarrow \infty$, and
\begin{eqnarray}
x_n(\delta) = x^\star(\delta)  + \frac{1}{\sqrt{n}} {W}_n(\delta) + \frac{1}{n}{V}_n(\delta)+ o_P(n^{-3/2}),
\label{eq:x1}
\end{eqnarray}
where ${W}_n(\delta)$ has mean $0$ and covariance matrix $ \xi(\delta) \equiv {\mathbb V}_{\mathbb P}[W_n(\delta)]$ and ${V}_n(\delta)$ is $O_P(1)$ with mean $\overline{V}(\delta)={\mathbb E}_{\mathbb P}[V_n(\delta)]$. $\xi(\delta)$ and $\overline{V}(\delta)$ are continuously differentiable in a neighborhood of $\delta=0$.
\end{proposition}
\begin{proof}
Uniqueness of $(x^\star(\delta), c^\star(\delta))$ was established in Proposition \eqref{prop:solution_differentiability}, while consistency of $(x_n(\delta), c_n(\delta))$ follows from  Theorem 5.9 in \cite{vdV}.
A proof of the remaining results can be found in the Appendix, where expressions for the random vectors $W_n(\delta)$ and $V_n(\delta)$, the covariance matrix $\xi(\delta)$ of $W_n(\delta)$ and the mean $\overline{V}(\delta)$ of $V_n(\delta)$ can also be found. 
\end{proof}

As in the case of the variance of $W_n(0)$ and the mean of $V_n(0)$ for SAA, $\xi(\delta)$ and $\overline{V}(\delta)$ do not depend on $n$. There is  a bias of $\frac{1}{n}\overline{V}(\delta)$ relative to the population DRO/DOO optimizer $x^\star(\delta)$ and the variance $\frac{1}{n}\xi(\delta)$ is different from that of the SAA solution.



\section{Out-of-sample expected reward} \label{sect:out-of-sample}

Let $x_n(\delta)$ be a solution of the best- or worst-case model, ${\mathbb E}_{\mathbb P}\big[f(x_n(\delta),\,Y_{n+1})\,|\,x_n(\delta)\big]$ the out-of-sample expected reward given $x_n(\delta)$, and
\begin{eqnarray*}
\mu_n(\delta) := {\mathbb E}_{\mathbb P}\big[f(x_n(\delta),\,Y_{n+1})\big] \equiv {\mathbb E}_{\mathbb P}\big\{{\mathbb E}_{\mathbb P}\big[f(x_n(\delta),\,Y_{n+1})\,|\,x_n(\delta)\big]\big\}
\end{eqnarray*}
be the expected out-of-sample reward after averaging over datasets $\{Y_1,\cdots,\,Y_n\}$ of size $n$ and outcome $Y_{n+1}$. Under the assumptions of our model, 
$Y_{n+1}$ has distribution $\mathbb P$, the distribution of $x_n(\delta)$ is characterized in Proposition \ref{prop:distribution-of-solution}, and $x_n(\delta)$ and $Y_{n+1}$ are independent. 
We now explore the relationship between $\mu_n(\delta)$ and the out-of-sample expected reward of the SAA optimizer, $\mu_n(0)$.

We denote the expected out-of-sample reward given $x$ by the function $g: {\mathbb R}^d \rightarrow {\mathbb R}$ defined by
\begin{eqnarray}
g(x) := {\mathbb E}_{\mathbb P}[f(x, Y)|x]. 
\label{eq:g}
\end{eqnarray}
Observe that $\mu_n(\delta) = {\mathbb E}_{\mathbb P}[g(x_n(\delta))]$. Concavity of $f(x, Y)$ in $x$ also implies that $g(x)$ is concave.

\subsection{Sample Average Approximation}\label{sec:SAA-stat}

Independence of $x_n(0)$ and $Y$, the concavity of $f(x, Y)$ and hence $g(x)$ in $x$, and Jensen's inequality imply that the out-of-sample expected reward under the SAA optimizer $x_n(0)$ is less than that of its mean ${\mathbb E}_{\mathbb P}[x_n(0)]$
\begin{eqnarray}
{\mathbb E}_{\mathbb P}[f(x_n(0), Y)] = {\mathbb E}_{\mathbb P}[g(x_n(0))] \leq g({\mathbb E}_{\mathbb P}[x_n(0)])= {\mathbb E}_{\mathbb P}\Big[f\big({\mathbb E}_{\mathbb P}[x_n(0)], Y\big)\Big] \leq g(x^\star(0)).
\label{eq:JG0}
\end{eqnarray}
This inequality is strict if $g(x) = {\mathbb E}_{\mathbb P}[f(x, Y)]$ is strictly concave  and $x_n(0)$ is random.


Equation \eqref{eq:JG0} suggests that we write the out-of-sample expected reward as a sum of the reward under the mean of the decision ${\mathbb E}_{\mathbb P}[x_n(0)]$ and a loss that we call the Jensen gap:
\begin{eqnarray}
\underbrace{{\mathbb E}_{\mathbb P}[g(x_n(0))]}_{\mu_n(0)}
=g({\mathbb E}_{\mathbb P}[x_n(0)]) + \underbrace{{\mathbb E}_{\mathbb P} \Big[g(x_n(0))- g({\mathbb E}_{\mathbb P}[x_n(0)])\Big]}_{\tiny \mbox{Jensen gap}}\leq  g(x^\star(0)).
\label{eq:decomp0}
\end{eqnarray}
The Jensen gap is negative when $g(x)$ is concave (and not linear) on the support of $x_n(0)$, while the optimality of $x^\star(0)$ for the population problem implies that ${\mathbb E}_{\mathbb P}[g(x_n(0))]$ is less than $g(x^\star(0))$ if  the solution is biased (${\mathbb E}_{\mathbb P}[x_n(0)] \neq x^\star(0)$).

Additional regularity (Assumptions \ref{ass:f} and \ref{ass:regf}) allows us to characterize the mean and variance of the SAA solution and to study  each component in the decomposition \eqref{eq:decomp0} separately.
Specifically,
we know from Proposition \ref{prop:SAA} that
\begin{eqnarray*}
x_n(0) = x^\star(0) + \frac{1}{\sqrt{n}}W_n(0) + \frac{1}{n}V_n(0) + o_P(n^{-1}).
\end{eqnarray*}
We make the assumption that the second-order error terms of the expectation and variance of $x_n(0)$ are $o(n^{-1})$. Sufficient conditions for this additional level of regularity are provided in Appendix \ref{app:Lp}.
\begin{assumption} \label{ass:SAALp}
The SAA solution is such that
\begin{eqnarray}
{\mathbb E}_{\mathbb P}[x_n(0)] & = & x^\star(0) + \frac{1}{n} \overline{V}(0) + o(n^{-1}), \label{eq:bias-var-SAA} \\
{\mathbb V}_{\mathbb P}[x_n(0)] & = & \frac{1}{n}\xi(0) + o(n^{-1}). \nonumber
\end{eqnarray}
\end{assumption}

Under Assumption \ref{ass:SAALp}, a Taylor series expansion of the first component of \eqref{eq:decomp0} around the population optimizer $x^\star(0)$ gives:
\begin{eqnarray}
g({\mathbb E}_{\mathbb P}[x_n(0)]) & = & g\Big(x^\star(0) + \frac{1}{n} \overline{V}(0) + o_P(n^{-1})\Big) 
\label{eq:decompTS-1} \\
& = & g(x^\star(0)) + \frac{1}{2n^2} \overline{V}(0)'\big[\nabla^2_x g(x^\star(0))\big]\overline{V}(0) + o(n^{-2})
\label{eq:decompTS1}
\end{eqnarray}
which shows that the
expected out-of-sample loss due to the finite-sample bias $\frac{1}{n} \overline{V}(0)$ is small, being on the order of $n^{-2}$. 
In the case of the Jensen gap, a Taylor series expansion around ${\mathbb E}_{\mathbb P}[x_n(0)]$ gives
\begin{eqnarray}
\lefteqn{{\mathbb E}_{\mathbb P} \Big[g(x_n(0))- g({\mathbb E}_{\mathbb P}[x_n(0)])\Big]} \nonumber \\
& = & \frac{1}{2}{\mathbb E}_{\mathbb P}\Big[\big(x_n(0) - {\mathbb E}_{\mathbb P}[x_n(0)]\big)'\nabla^2_x g({\mathbb E}_{\mathbb P}[x_n(0)])\big(x_n(0) - {\mathbb E}_{\mathbb P}[x_n(0)]\big)\Big] +o(n^{-1}) \nonumber \\
& = & \frac{1}{2}{\rm tr}\Big\{{\mathbb E}_{\mathbb P}\Big[\big(x_n(0) - {\mathbb E}_{\mathbb P}[x_n(0)]\big)\big(x_n(0) - {\mathbb E}_{\mathbb P}[x_n(0)]\big)'\Big] \nabla^2_x g\Big(x^\star(0) + \frac{1}{n} \overline{V}(0) + o(n^{-1})\Big) \Big\}+o(n^{-1}) \nonumber \\
& = & \frac{1}{2n}\mbox{tr} \big\{\xi(0) \nabla^2_x g(x^\star(0))\big\} + o(n^{-1}),
\label{eq:decompTS2}
\end{eqnarray}
where the second and third equalities follow from \eqref{eq:bias-var-SAA} and \eqref{eq:decompTS-1}.
The Jensen gap is negative because $g(x)$ is concave
and becomes more negative when the variance $\frac{1}{n}\xi(0)$ of $x_n(0)$ increases or the curvature $\nabla^2_x g(x^\star(0))$ becomes more negative.

Adding  \eqref{eq:decompTS1} and \eqref{eq:decompTS2} we obtain an expression for the out-of-sample expected reward in terms of  the population optimizer:
\begin{eqnarray}
\underbrace{{\mathbb E}_{\mathbb P}[f(x_n(0), Y)] }_{\tiny \mu_n(0)} &= & {\mathbb E}_{\mathbb P}[f(x^\star(0), Y)]  + \frac{1}{2n}\underbrace{\mbox{tr}\big\{\xi(0){\mathbb E}_{\mathbb P}[\nabla^2_xf(x^\star(0), Y)]\big\}}_{\tiny {\mathbb E}_{\mathbb P}[W_n(0)' \nabla^2_xf(x^\star(0), Y) W_n(0)]<0}   +
o(n^{-1}).
\label{eq:SAA_reward}
\end{eqnarray}
This expression shows that the loss in expected reward is of order $n^{-1}$ and comes from the uncertainty in the solution through the Jensen gap. There is also a loss due to the finite-sample bias in the SAA solution, but this is an order of magnitude smaller.



\subsection{Distributionally Robust/Optimistic Optimization} \label{sec:DRO}

To compare the out-of-sample expected reward of the DRO solution with that of SAA, we write ${\mathbb E}_{\mathbb P}[g(x_n(\delta))]$ in terms of the SAA reward ${\mathbb E}_{\mathbb P}[g(x_n(0))]$ and additional terms that show how changes in the mean and variance of the solution affect the expected reward:
\begin{eqnarray}
{\mathbb E}_{\mathbb P}[g(x_n(\delta))] & = & {\mathbb E}_{\mathbb P}[g(x_n(0))] - \overbrace{\Big\{g({\mathbb E}_{\mathbb P}[x_n(0)]) - g({\mathbb E}_{\mathbb P}[x_n(\delta)])\Big\} }^{\tiny \mbox{(A) Loss from change in mean of solution}}\nonumber \\ [5pt]
& & - \Big\{\underbrace{{\mathbb E}_{\mathbb P}\Big[g(x_n(0)) - g({\mathbb E}_{\mathbb P}[x_n(0)])\Big] - {\mathbb E}_{\mathbb P}\Big[g(x_n(\delta)) - g({\mathbb E}_{\mathbb P}[x_n(\delta)])\Big]}_{\tiny \mbox{(B) Change in the Jensen gap}}\Big\}.
\label{eq:DRO-decomp}
\end{eqnarray}
The term (A) shows how a change in the expected value of the solution changes the reward,  while (B) is the associated change in the Jensen gap.

 DRO and DOO add a bias and change the variance of the SAA optimizer.
Specifically, it follows from Proposition \ref{prop:distribution-of-solution} and equations \eqref{eq:pop_rob_asymp_bias} and \eqref{eq:bias-var-SAA} (under Assumptions \ref{ass:f},  \ref{ass:phi} that \eqref{eq:x1} holds. As in the case of the SAA solution, we make the additional assumption  that the second-order error terms of the mean and variance of $x_n(\delta)$ is $o(n^{-1})$. A sufficient condition  is given in Appendix \ref{app:Lp}.
\begin{assumption} \label{ass:DLp}
There exists a neighborhood of $\delta=0$ such that
\begin{eqnarray}
 {\mathbb E}_{\mathbb P}[x_n(\delta)] & = & x^\star(\delta)+\frac{1}{n}\overline{V}(\delta)  + o(n^{-1})
 \label{eq:bias-var-DRO}
  \\ [5pt]
 {\mathbb V}_{\mathbb P} [x_n(\delta)] & = & \frac{1}{n}\xi(\delta) + o(n^{-1}).
 \nonumber
\end{eqnarray}
\end{assumption}
Under this assumption and for $n$ sufficiently large, there is an open neighborhood of $\delta=0$ such that
\begin{eqnarray}
{\mathbb E}_{\mathbb P}[x_n(\delta)] & = & {\mathbb E}_{\mathbb P}[x_n(0)] + \delta \Big(\pi + \frac{1}{n}\overline{V}_\delta(0)\Big) + O(\delta^2) + \frac{1}{n} O(\delta^2) + o(n^{-1})
\label{eq:mean-var}  \\
{\mathbb V}_{\mathbb P} [x_n(\delta)] & = & {\mathbb V}_{\mathbb P}[x_n(0)] + \frac{\delta}{n}\xi_\delta(0) + \frac{1}{n} O(\delta^2) + o(n^{-1}). \nonumber
\end{eqnarray}
Here, $\overline{V}_\delta(0) = \frac{\mathrm{d}}{\mathrm{d}\delta}\overline{V}(0)$ and $\xi_\delta(0) = \frac{\mathrm{d}}{\mathrm{d}\delta}\xi(0)$ denote derivatives of $\xi(\delta)$ and $\overline{V}(\delta)$ at $\delta=0$, which exist because $\xi(\delta)$ and $\overline{V}(\delta)$ are continuously differentiable at $\delta=0$  (Proposition \ref{prop:distribution-of-solution}).

To the first order, the change in the bias relative to the SAA solution is $\delta\big(\pi + \frac{1}{n}\overline{V}_\delta(0)\big)$ and $\frac{\delta}{n} \xi_\delta(0)$ is the change in the variance. The expression for ${\mathbb E}_{\mathbb P}[x_n(\delta)]$ follows from the previous equation together with \eqref{eq:pop_rob_asymp_bias} and \eqref{eq:bias-var-SAA} and shows that  the difference between the mean of $x_n(\delta)$ and $x_n(0)$ comes from two sources, the change in the population solution from $x^\star(0)$ to $x^\star(\delta)$ due to robustness which contributes the $\pi$, and the change in the finite-sample bias which gives $\frac{1}{n}\overline{V}_\delta(0)$. The variance of the solution also changes by $\frac{1}{n}\xi_\delta(0)$. An explicit relationship between terms (A) and (B) in \eqref{eq:DRO-decomp} and the robustness parameter can be  obtained by substituting  \eqref{eq:mean-var} into each term in \eqref{eq:DRO-decomp}.

For the first term (A), it is shown in Appendix \ref{App:DRO-obj} that
\begin{eqnarray}
\lefteqn{g({\mathbb E}_{\mathbb P}[x_n(\delta)]) - g({\mathbb E}_{\mathbb P}[x_n(0)]) } \label{eq:DRO-exp1} \\ [5pt]
& =&   \underbrace{\frac{\delta^2}{2}\pi'\big[\nabla^2_x g(x^\star(0))\big]\pi  + \frac{\delta}{n} \Big\{ \pi'\nabla^2_x g(x^\star(0))\overline{V}(0)\Big\}}_{\tiny \mbox{Impact of robustness}}
 + o(1/n^2)+ o(\delta^2) + \frac{1}{n} o(\delta^2) + \frac{1}{n^2}O(\delta). \nonumber
\end{eqnarray}

In the case of the Jensen gap (B):
\begin{eqnarray}
\lefteqn{{\mathbb E}_{\mathbb P} \Big\{g(x_n(\delta))- g({\mathbb E}_{\mathbb P}[x_n(\delta)])\Big\}- {\mathbb E}_{\mathbb P} \Big\{g(x_n(0))- g({\mathbb E}_{\mathbb P}[x_n(0)])\Big\} }  \label{eq:DRO-exp2} \\ [5pt]
& = & \underbrace{\Big(\frac{\delta}{2n} \Big)\frac{\mbox{d}}{\mbox{d}\delta}\Big\{\mbox{tr}\Big(\nabla^2_x g(x^\star(\delta))\xi(\delta)\Big)\Big\}\Big|_{\delta=0}}_{\tiny \mbox{Impact of robustness}}+ \frac{1}{n}O(\delta^2) + O(n^{-3/2})
\nonumber
\end{eqnarray}
where
\begin{eqnarray*}
\frac{\mbox{d}}{\mbox{d}\delta}\Big\{\mbox{tr}\Big(\nabla^2_x g(x^\star(\delta))\xi(\delta)\Big)\Big\}\Big|_{\delta=0} = \pi' \nabla_x\Big(\mbox{tr}\Big\{\nabla^2_x g(x)\xi(0)\Big\}\Big)\Big|_{x=x^\star(0)} + \mbox{tr}\Big(\nabla^2_xg(x^\star(0))\,\xi_\delta(0)\Big).
\end{eqnarray*}
This expression shows that the change in the Jensen gap comes from the change in the curvature of the objective function due to the DRO/DOO bias (first term) and the change in the variance of the solution (second term).




The out-of-sample expected reward under a DRO/DOO solution can be written in terms of the expected reward under the SAA optimizer by adding \eqref{eq:DRO-exp1} and \eqref{eq:DRO-exp2}.



\begin{proposition}\label{prop:out-of-sample reward}
Suppose $\{Y_1,\cdots,\,Y_n\}$ and $Y_{n+1}$ are drawn iid from $\mathbb P$,  $f(x, Y)$ satisfies Assumption \ref{ass:f} and \ref{ass:regf}, $\phi(z)$ satisfies Assumption \eqref{ass:phi}, and
Assumptions \ref{ass:reg}, \ref{ass:SAALp} and \ref{ass:DLp} hold.
Then for every $n$ sufficiently large, there is an open neighborhood of $\delta=0$ such that
\begin{eqnarray}
\lefteqn{{\mathbb E}_{\mathbb P}\big[f(x_n(\delta),\,Y_{n+1})\big]} \nonumber \\
& = & {\mathbb E}_{\mathbb P}\big[f(x_n(0),\,Y_{n+1})\big] +  \delta \frac{\rho}{n}  + \frac{\delta^2}{2} \pi' {\mathbb E}_{\mathbb P}\big[\nabla_x^2 f(x^{\star}(0),\,Y_{n+1})\big]\pi + o(\delta^2) + o(n^{-1})
\label{eq:out-of-sample-reward}
\end{eqnarray}
where 
\begin{eqnarray}
\rho & = & \overline{V}(0)'\,{\mathbb E}_{\mathbb P}[\nabla^2_xf(x^\star(0), Y)]\,\pi + \frac{1}{2}\, \frac{\mathrm{d}}{\mathrm{d}\delta}\mathrm{tr}\Big(\xi(\delta){\mathbb E}_{\mathbb P}[\nabla^2_xf(x^\star(\delta), Y)]\Big)\Big|_{\delta=0} \nonumber \\[5pt]
&= & \overline{V}(0)'\,{\mathbb E}_{\mathbb P}[\nabla^2_xf(x^\star(0), Y)]\,\pi \nonumber \\
& & + \frac{1}{2}\mathrm{tr}\Big(\xi_\delta(0){\mathbb{E}_{{\mathbb{P}}}[\nabla^2_x f(x^\star(0),\,Y)]}\Big) +\frac{1}{2} \pi'\nabla_x \Big[\mathrm{tr}\Big(\xi(0)\,{\mathbb E}_{{\mathbb P}}\big[\nabla_x^2 f(x^{\star}(0),\,Y_{n+1})\big]\Big)\Big] \label{eq:rho}
\end{eqnarray}
and $\pi$ is given by \eqref{eq:pop_pi}.
\end{proposition}

It is worth reiterating the origins of each term in $\rho$. The first is the change in the expected reward \eqref{eq:DRO-exp1} that results from the change to the mean of the SAA solution, while the second is the change in the Jensen gap that results from the change in the variance of the solution. The last term is the change in the Jensen gap resulting from a change in the curvature of the objective function due to the DRO/DOO bias \eqref{eq:DRO-exp2}. When $f(x, Y)$ is concave but not linear, it is clear from \eqref{eq:rho} that $\rho$ is unlikely to be zero, except in pathological cases.



The following result shows that it is always possible to select $\delta$ so that the out-of-sample expected reward under the DRO/DOO solution $x_n(\delta)$ exceeds that of $x_n(0)$ and provides an estimate of the size of the outperformance.
\begin{theorem}\label{theorem:out-of-sample}
Suppose that the assumptions of Proposition \ref{prop:out-of-sample reward} hold.
If $\rho\neq 0$, then for $n$ sufficiently large, $\delta$ can always be selected so that the expected reward under $x_n(\delta)$ exceeds that of  the SAA  optimizer out-of-sample
\begin{eqnarray*}
{\mathbb E}_{\mathbb P}\big[f(x_n(\delta),\,Y_{n+1})\big]>{\mathbb E}_{\mathbb P}\big[f(x_n(0),\,Y_{n+1})\big].
\end{eqnarray*}
If $n$ is such that the constant
\begin{eqnarray*}
\delta_n = -\frac{1}{n}\frac{\rho}{\pi'{\mathbb E}_{\mathbb P}[\nabla^2_xf(x^\star(0), Y)] \pi}
\end{eqnarray*}
is in the open neighborhood of $\delta=0$ where \eqref{eq:out-of-sample-reward} holds, then
\begin{eqnarray}
\lefteqn{{\mathbb E}_{\mathbb P}\big[f(x_n(\delta_n),\,Y_{n+1})\big] - {\mathbb E}_{\mathbb P}\big[f(x_n(0),\,Y_{n+1})\big]} \nonumber \\
& = &  - \frac{1}{2n^2} \underbrace{\Big(\frac{\rho^2}{\pi'{\mathbb E}_{\mathbb P}[\nabla^2_xf(x^\star(0), Y)] \pi}\Big)}_{\tiny \mbox{negative}} + o(n^{-2})
\label{eq:opt-outperform} \\
& > &  0 \nonumber
\end{eqnarray}
\end{theorem}

\begin{proof}
Equation \eqref{eq:out-of-sample-reward} shows how the out-of-sample expected reward depends on the robustness parameter $\delta$ when it is small. If $\rho$ is non-zero,  the change in the expected reward
\begin{eqnarray*}
\lefteqn{{\mathbb E}_{\mathbb P}\big[f(x_n(\delta),\,Y_{n+1})\big] - {\mathbb E}_{\mathbb P}\big[f(x_n(0),\,Y_{n+1})\big]} \\ & = & \delta \frac{\rho}{n}  + \frac{\delta^2}{2} \pi' {\mathbb E}_{\mathbb P}\big[\nabla_x^2 f(x^{\star}(0),\,Y_{n+1})\big]\pi + o(\delta^2) + o(1/n)
\end{eqnarray*}
is dominated by the linear term $\delta \frac{\rho}{n}$ when $\delta$ is small and we can always choose $\delta$ so that  the expected reward from $x_n(\delta)$ exceeds that from the SAA solution $x_n(0)$. If $\rho$ is positive, a positive choice of $\delta$  (DRO) does the job, while a negative value of $\delta$ (DOO) should be chosen if $\rho$ is negative. The expression for $\delta_n$ is obtained by maximizing the quadratic on the right-hand-side of \eqref{eq:out-of-sample-reward}. 

The inequality \eqref{eq:opt-outperform}  follows by substituting the expression for $\delta_n$ into \eqref{eq:out-of-sample-reward}.
\end{proof}

Although we can always out-perform SAA,  \eqref{eq:opt-outperform} suggests that it is on the order of $n^{-2}$, regardless of whether it was achieved by best-case or worst-case optimization. Out-performance from either DRO or DOO is unlikely to be large if $n$ is large or $\rho$ is small. More generally, $\frac{\rho}{n}$ is the sensitivity of the sum of (A) and (B) in \eqref{eq:DRO-decomp} to changes in $\delta$. Our assumptions allow us to estimate this explicitly. (A) and (B) may depend differently on $n$ and $\delta$ in applications where our assumptions do not hold, though we have not explored this issue.

It is difficult to compute $\rho$  because it depends on the solution of the population problem $x^\star(0)$ and the population distribution $\mathbb P$, both of which are unknown, and hence to determine {\it a priori} whether the optimistic or worst-case solution will out-perform SAA. This can be sidestepped by optimizing a bootstrap or cross-validation estimate of the out-of-sample reward over $\delta$, which we explore in an inventory problem in Section \ref{sect:wcs}.

One limitation of our analysis is that we make strong assumptions about the differentiability of $f(x, Y)$. These assumptions enable us to characterize the statistical properties of the SAA and DRO/DOO solutions and to analyze out-of-sample performance. However, there are many applications where they are not satisfied, the inventory problem below being one such case. Our assumptions enable us to illustrate what is possible with DOO and to provide insight into why it occurs. 
We do not rule out the possibility that either DRO or DOO will out-perform SAA when the differentiability assumptions we impose are not satisfied.
All predictions that stem from our analysis, both in Theorem \ref{theorem:out-of-sample} and later sections, are observed in the inventory example, where the objective function violates Assumption \ref{ass:f}.


\begin{example} \label{ex:AP}
In \cite{AP}, many examples where the out-of-sample expected reward from DRO exceeds that of  SAA are presented for an objective function of the form
\begin{eqnarray*}
f(x, Y) = \frac{1}{2}  x' H x + v(Y)'x + u(Y),
\end{eqnarray*}
where $H$ is strictly negative definite. Since ${\mathbb E}_{\mathbb P}[\nabla^2_xf(x^{\star}(0),\, Y)]=H$, ${\mathbb E}_{\mathbb P}[\nabla_xf(x^{\star}(0),\, Y)]=0$ and all derivatives of $f(x, Y)$ of order $3$ or higher are zero, it follows from \eqref{eq:V-bar} that $\overline{V}(0) = 0$ and
\begin{eqnarray*}
\frac{1}{2} \pi'\nabla_x \Big[\mathrm{tr}\Big(\xi(0)\,{\mathbb E}_{{\mathbb P}}\big[\nabla_x^2 f(x^{\star}(0),\,Y_{n+1})\big]\Big)\Big] = 0,
\end{eqnarray*}
and hence by \eqref{eq:rho}
\begin{eqnarray*}
\rho = \frac{1}{2}\mathrm{tr}\Big(\xi_\delta(0){\mathbb{E}_{{\mathbb{P}}}[\nabla^2_x f(x^\star(0),\,Y)]}\Big) = \frac{1}{2}{\rm tr}\Big(\xi_\delta(0)H\Big).
\end{eqnarray*}
The only term from $\rho$ that remains corresponds to the change in the Jensen gap due to a change in the variance of the solution.

If $x$ is one-dimensional,  $H$ is a negative constant.
If DRO reduces the variance of the solution relative to SAA, $\xi_\delta(0)$ is negative and $\rho=\frac{1}{2}\xi_\delta(0)H$ is positive, so  DRO will also out-perform SAA if $\delta>0$ is chosen appropriately.
If DRO increases the variance of the solution (see \cite{GKL-cal} for an example where this happens), $\rho$ is negative and a DOO solution will have a higher out-of-sample expected reward than SAA.

 All examples in \cite{AP} where DRO out-performs SAA are ones where the change in variance is such that $\rho=\frac{1}{2}{\rm tr}(\xi_\delta(0)H)$ is positive. Optimistic optimization will beat SAA in the examples from \cite{AP} when DRO does not.

\end{example}


\begin{example} \label{ex:inv}
Consider an inventory problem with reward
\begin{eqnarray}
f(x, Y)= r \min\{x,  Y\} +  q\max\{x-Y, 0\}- s\max\{Y-x, 0\} - c x
\label{eq:inv2}
\end{eqnarray}
where the selling price $r=10$, the purchase cost $c=9$, and the shortage and scrap costs, $s$ and $q$, are zero. Under the population distribution, demand $Y = \max\{m + I X_1 - (1-I) X_2, 0\}$, where $X_i$ is exponential with mean $\mu_i$ ($i=1, 2$), $I$ is Bernoulli with $p=P[I=1]$, and $m$ is a constant; $X_1, X_2$ and $I$ are independent. In this experiment, $m=250$, $\mu_1= 10$, $\mu_2=60$ and $p=0.9$.

We note that the objective function \eqref{eq:inv2} does not satisfy the differentiability assumptions that facilitate our analysis  (Assumptions \ref{ass:f}, \ref{ass:regf}, and \ref{ass:reg}). It will be seen, however, that all predictions from our analysis will be observed in this inventory model. More generally, this provides some evidence that the insights from our analysis are true under less stringent assumptions about the degree of smoothness of the objective function.


We generated $5,000$ datasets, each of size $n=30$, and solved the DRO/DOO problems with a modified-$\chi^2$  penalty function ($\phi(z) = \frac{1}{2}(z-1)^2$) over a range of $\delta$ for each dataset. We then computed the out-of-sample expected reward for the family of DRO/DOO decisions $x_n(\delta)$, giving $5,000$ conditional expected-reward functions. The out-of-sample expected reward $\mu_n(\delta) = {\mathbb E}_{\mathbb P}[f(x_n(\delta), Y_{n+1})]$ shown in Figure \ref{fig:mean} was estimated by averaging these $5,000$ functions.

If we restrict ourselves to worst-case solutions (DRO), the SAA solution ($\delta=0$) is optimal and the out-of-sample expected reward is $189$.
However, the out-of-sample expected reward is maximized with an optimistic solution $x_n(\delta)$ with $\delta=-1.3\times 10^{-3}$ and has value $193$.
With different parameter values and/or population distribution, a worst-case model $(\delta>0$)  could be optimal.

\begin{figure}[h]
\includegraphics[scale=0.25]{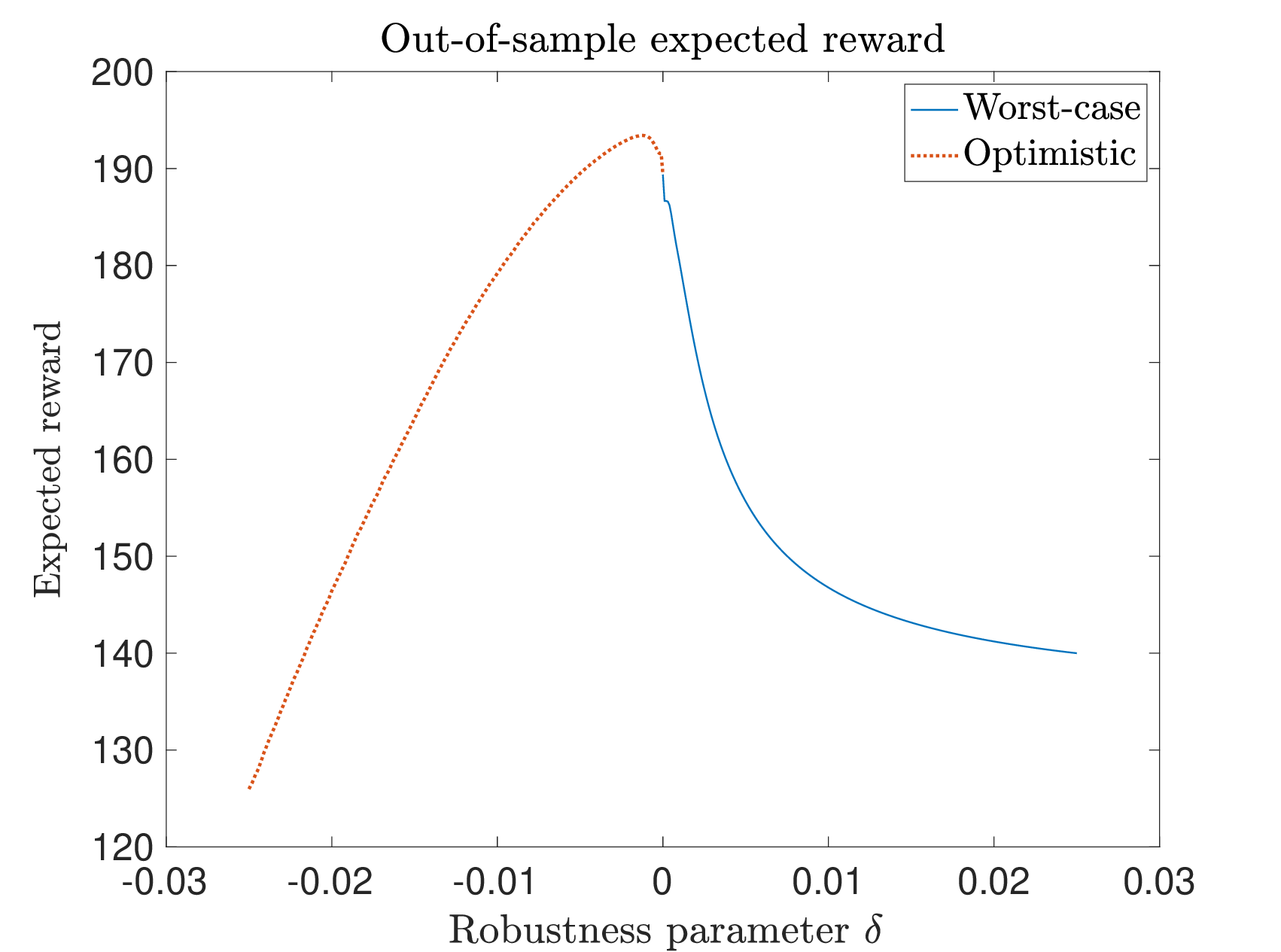}
\caption{The plot shows the out-of-sample expected reward $\mu_n(\delta) = {\mathbb E}_{\mathbb P}[f(x_n(\delta), Y_{n+1})]$ as a function of $\delta$. The optimistic solution $x_n(\delta)$ with $\delta=-1.3\times 10^{-3}$ maximizes the expected reward  ($193$). The expected reward associated with SAA  ($\delta=0$) is $189$. }
\label{fig:mean}
\end{figure}

\end{example}

In summary, DRO models typically come with a free parameter -- the size of the uncertainty set or the robustness parameter $\delta$ -- which many select by optimizing an estimate of the out-of-sample expected reward via cross-validation or the bootstrap. However, it is not always possible to beat SAA if we restrict ourselves to worst-case models \cite{AP,GKL-cal}. If the ultimate goal is to beat SAA out-of-sample, it is reasonable to consider worst- and best-case models to ensure  this is  possible.

\section{So, what is the catch?} \label{sect:wcs}

\subsection{Best-case solutions are not robust} \label{sec:WCS}
We have focused on the out-of-sample expected reward of solutions of best- and worst-case problems. 
However,
it has also been shown  that
DRO is intrinsically a tradeoff between (in-sample) expected reward and worst-case sensitivity \cite{GKL-cal,GKL,GKL-WCS}.  Worst-case sensitivity is a quantitative measure of robustness and it  is natural to evaluate the worst-case sensitivity of the best-case optimizer.

We recall the definition of worst-case sensitivity \cite{GKL}. Suppose decision $x$ is fixed and ${\mathbb E}_{{\mathbb P}_n}[f(x, Y)]$ be the expected reward under the nominal distribution ${\mathbb P}_n$. Worst-case sensitivity ${\mathcal S}_{{\mathbb P}_n}(f(x,\,Y))$ is the rate of decrease of the in-sample expected reward under ``worst-case deviations" from the nominal distribution.  When the difference between probability distributions is measured by smooth $\phi$-divergence, it is natural to define the worst-case distribution as the minimizer
\begin{align*}
{\mathbb Q}(\varepsilon) &
:= \left\{
\begin{array}{cl}
\begin{displaystyle}\arg\min_{\mathbb Q}\Big\{  \sum_{i=1}^nq_i f(x,\,Y) +  \frac{1}{\varepsilon} \sum_{i=1}^n {p}^n_i \phi\Big(\frac{q_i}{{p}^n_i}\Big)\Big\}, \end{displaystyle}& \varepsilon>0,\\ [10pt]
{\mathbb P}_n, &\varepsilon=0.
\end{array}\right.
\end{align*}
Here, the penalty $\sfrac{1}{\varepsilon}>0$ on $\phi$-divergence controls the size of the deviation of ${\mathbb Q}(\varepsilon)$ from ${\mathbb P}_n$, which is increasing in $\varepsilon$. Worst-case sensitivity is
\begin{eqnarray}
{\mathcal S}_{{\mathbb P}_n}(f(x,\,Y)) & :=  & - \frac{\mathrm{d}}{\mathrm{d}\varepsilon}{\mathbb E}_{{\mathbb Q}(\varepsilon)}[f(x,\,Y)]\Big|_{\varepsilon=0} \nonumber  \\[5pt]
& = & - \lim_{\varepsilon\downarrow 0}\frac{{\mathbb E}_{{\mathbb Q}(\varepsilon)}[f(x,\,Y)] -{\mathbb E}_{{\mathbb P}_n}[f(x,\,Y)]}{\varepsilon}  \nonumber  \\[5pt]
& = & \frac{1}{\phi''(1)}{\mathbb V}_{{\mathbb P}_n}[f(x,\,Y)], \label{eq:wcs}
\end{eqnarray}
which, in the case of smooth $\phi$-divergence, is equal to the in-sample variance of the reward. There are other ways to define worst-case sensitivity: we can use a different divergence measure, or we can control the ``distance" of $\mathbb Q$ from ${\mathbb P}_n$ through a constraint and look at the limit when it vanishes \cite{GKL-WCS}. However, they all capture a similar idea.

Worst-case sensitivity is an in-sample {\it measure of robustness}. Given a nominal model ${\mathbb P}_n$ and a decision $x$, it quantifies the sensitivity of the expected reward to errors in the nominal model and is one way to evaluate whether a given decision is more or less robust than another.
Worst-case sensitivity depends on the choice of uncertainty set, though it is always a {\it generalized measure of deviation (spread)} \cite{GKL-WCS}. Intuitively, the expected reward is sensitive to  mis-specification when it has a large spread because small changes in the probability of extreme rewards  (positive {\it or} negative) can have a big impact on the mean.

Instead of $x$, we can substitute the solution $x_n(\delta)$ of the best/worst-case problem into \eqref{eq:wcs}. To see the impact of best/worst case optimization on the ``robustness" of the solution, we write the sensitivity of $x_n(\delta)$ in terms of the sensitivity of the SAA solution $x_n(0)$.
Using the expansion \eqref{eq:insample_rob_asymp_bias}, the in-sample variance of the reward under  $x_n(\delta)$ is
\begin{eqnarray*}
{\mathbb V}_{{\mathbb P}_n}[f(x_n(\delta),\,Y)]
=  {\mathbb V}_{{\mathbb P}_n}[f(x_n(0),\,Y)]  + 2 \frac{\delta}{\phi^{''}(1)}\, \beta_n'\Big({\mathbb{E}}_{{\mathbb P}_n}\big[\nabla_x^2 f(x_n(0),\,Y) \big]\Big)^{-1}\beta_n + O(\delta^2)
\end{eqnarray*}
where
\begin{eqnarray*}
\beta_n & = &  \mathrm{Cov}_{\mathbb{P}_n}\Big[{\nabla_x f(x_n(0),\,Y)}, \, f(x_n(0),\,Y)\Big].
\end{eqnarray*}
It now  follows from \eqref{eq:wcs} that worst-case sensitivity
\begin{eqnarray}
\lefteqn{{\mathcal S}_{{\mathbb P}_n}(f(x_n(\delta),\,Y))}\nonumber  \\
 &  = &  \underbrace{{\mathcal S}_{{\mathbb P}_n}(f(x_n(0),\,Y))}_{\tiny \mbox{Sensitivity of SAA solution}}+ 2 \frac{\delta}{[\phi^{''}(1)]^2}\, \beta_n'\Big({\mathbb{E}}_{{\mathbb P}_n}\big[\nabla_x^2 f(x_n(0),\,Y) \big]\Big)^{-1}\beta_n + O(\delta^2).
\label{eq:sensitivity}
\end{eqnarray}
Strict concavity of $f(x, Y)$ in $x$ implies that
\begin{eqnarray*}
\beta_n'\Big({\mathbb{E}}_{{\mathbb P}_n}\big[\nabla_x^2 f(x_n(0),\,Y) \big]\Big)^{-1}\beta_n < 0
\end{eqnarray*}
so solutions of the optimistic DOO problem ($\delta<0$) have a larger sensitivity than SAA, while worst-case solutions ($\delta>0$) have a smaller sensitivity and are therefore more robust.

\subsection{The optimal value of $\delta$ may be difficult to estimate} \label{sec:est-delta}

Although there exists an optimal DOO or DRO decision that has a larger out-of-sample expected reward than SAA, the parameter $\delta$ corresponding to this decision  needs to be estimated. It is natural to use bootstrap or cross-validation, though the estimate will depend on the data set so there will be estimation error. If the optimal value of $\delta$ is positive but we erroneously select a negative value of $\delta$ (i.e. DOO), not only do we lose a little out-of-sample expected reward but we increase worst-case sensitivity and the solution will be less robust than SAA.

\begin{example}[Inventory control (continued)]
For each of the $5,000$ sampled datasets (each with $n=30$ datapoints), we compute worst-case sensitivity ${\mathcal S}_{{\mathbb P}_n}(f(x_n(\delta),\,Y))$ as a function of $\delta$. Figure \ref{fig:sensitivityA} is obtained by averaging these. As shown in \eqref{eq:sensitivity} optimistic solutions have a larger sensitivity than SAA and DRO solutions in the neighborhood of $\delta=0$ and sensitivity is linear in $\delta$.  For each of the $5,000$ datasets, we can compute the in-sample mean-sensitivity frontier corresponding to the DRO/DOO solutions over a range of $\delta$.  Figure \ref{fig:sensitivityB} is obtained by averaging these frontiers.

Figure \ref{fig:oos_frontier} shows that out-of-sample mean-variance frontier. SAA can be beaten by an optimistic decision, but this comes at the cost of a large increase in the out-of-sample sensitivity (variance).

Figure \ref{fig:delta_plots} shows the distribution of bootstrap estimates of the optimal $\delta$ as a function of the number of data points. Here we simulated $500$ independent datasets of size $n=15$ and $n= 30$, and used 50 boostrap samples for each dataset to estimate of the out-of-sample expected reward which we optimized over $\delta$. When $n=15$, $\delta=-1.6\times 10^{-3}$ is optimal; when $n=30$, $\delta=-1.3\times10^{-3}$.
The bootstrap estimates are not very accurate when $n=15$ (mean = $1.4\times10^{-2}$,  SD = $2.8\times10^{-2}$), with only $49\%$ of experiments giving an estimate with the correct sign, and there is a long right tail. This is not very surprising given that a sample of size $n=15$ is  small when the population demand has a standard deviation of $28$. The chance of getting a data set that is not representative of the population is high, leading to poor estimates of $\delta$ (e.g. the extreme values in the right tail of the distribution). Accuracy improves with $n=30$ data points (mean = $-1.3\time10^{-2}$, SD = $4.7\times10^{-3}$), with estimates having the correct sign $88\%$ of the time and a smaller skew. It could well be the case that other methods estimate the optimal $\delta$, or at least its sign, more accurately. We have not explored this issue further.



\begin{figure}[ht]
\centering
\begin{subfigure}{.495\textwidth}
\centering
\includegraphics[scale=0.225]{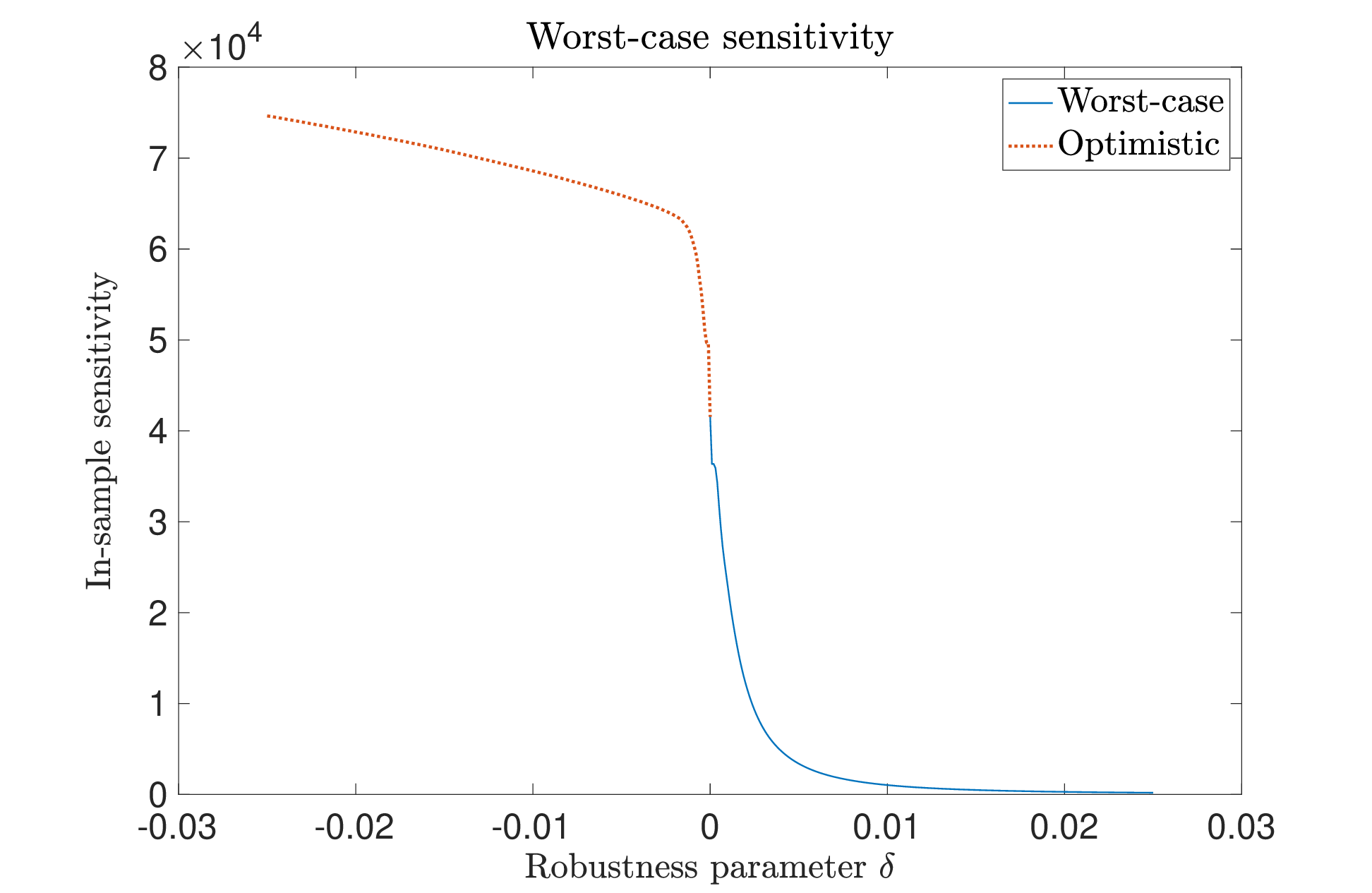}
\caption{Worst-case sensitivity}
\label{fig:sensitivityA}
\end{subfigure}
\begin{subfigure}{.495\textwidth}
\centering
\includegraphics[scale=0.225]{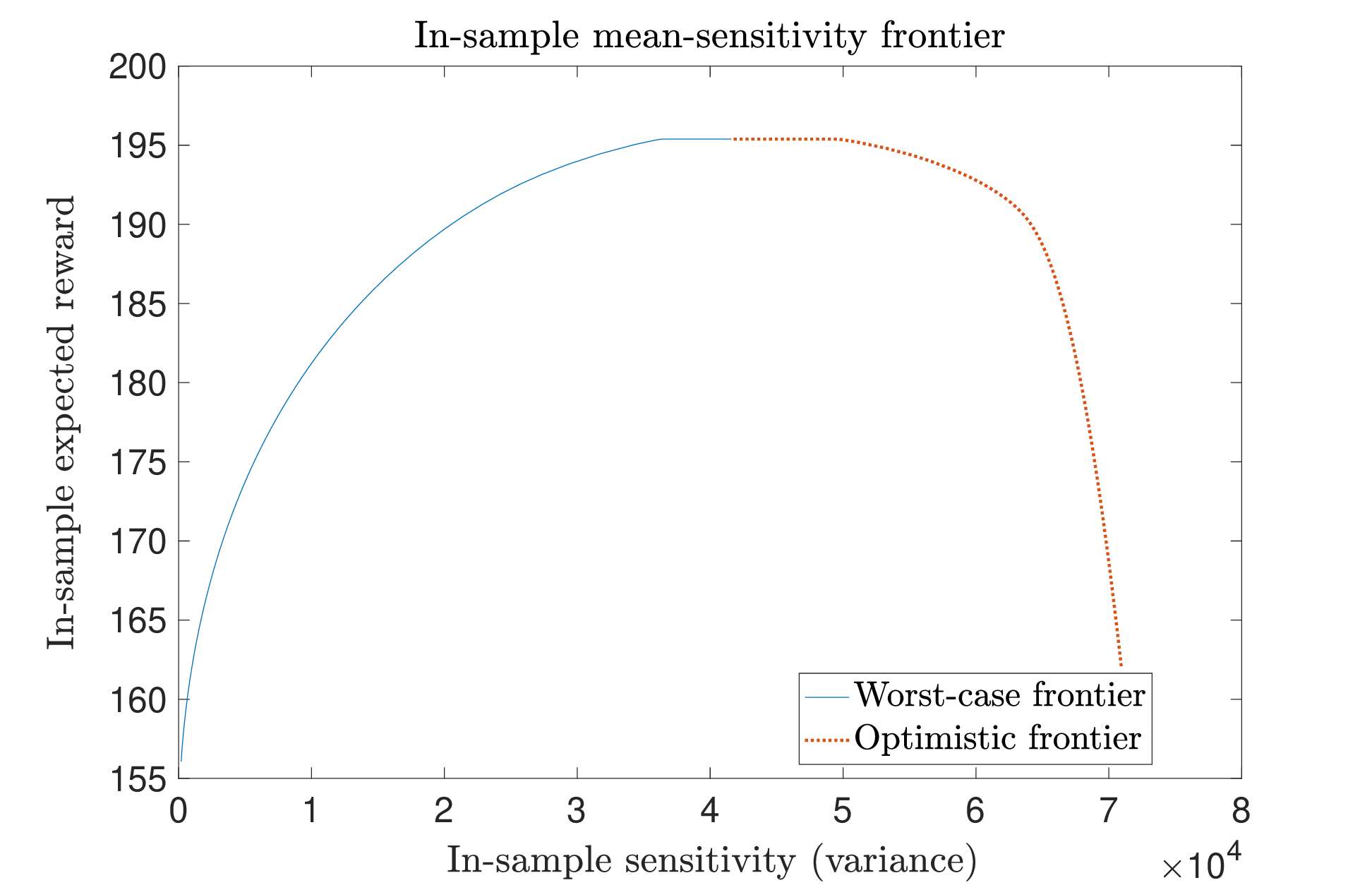}
\caption{Average in-sample mean-sensitivity frontier}
\label{fig:sensitivityB}
\end{subfigure}

\vspace{0.75cm}

\begin{subfigure}{.495\textwidth}
\centering
\includegraphics[scale=0.225]{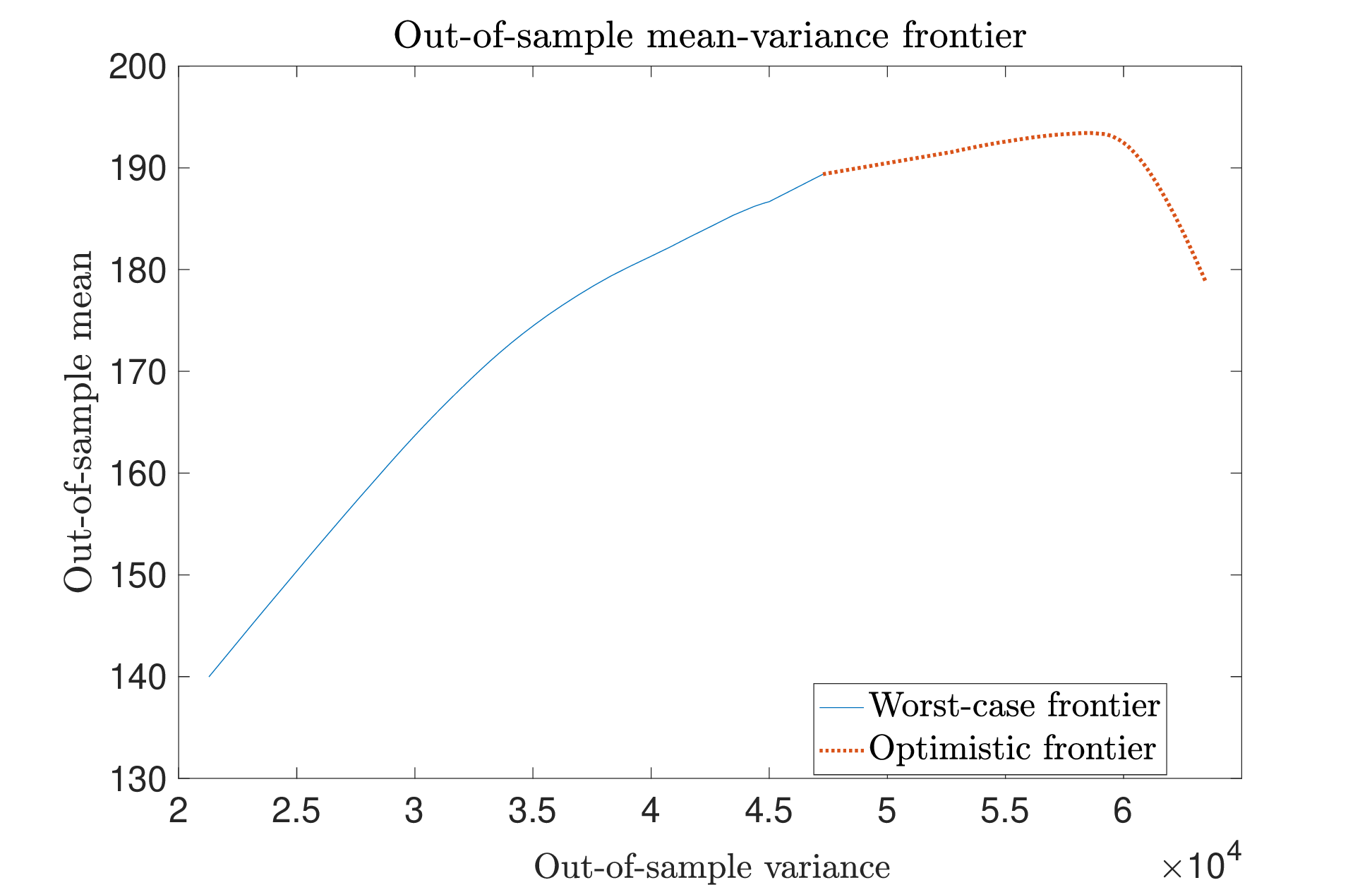}
\caption{Out-of-sample mean-sensitivity frontier.}
\label{fig:oos_frontier}
\end{subfigure}
\begin{subfigure}{.495\textwidth}
\centering
\includegraphics[scale=0.225]{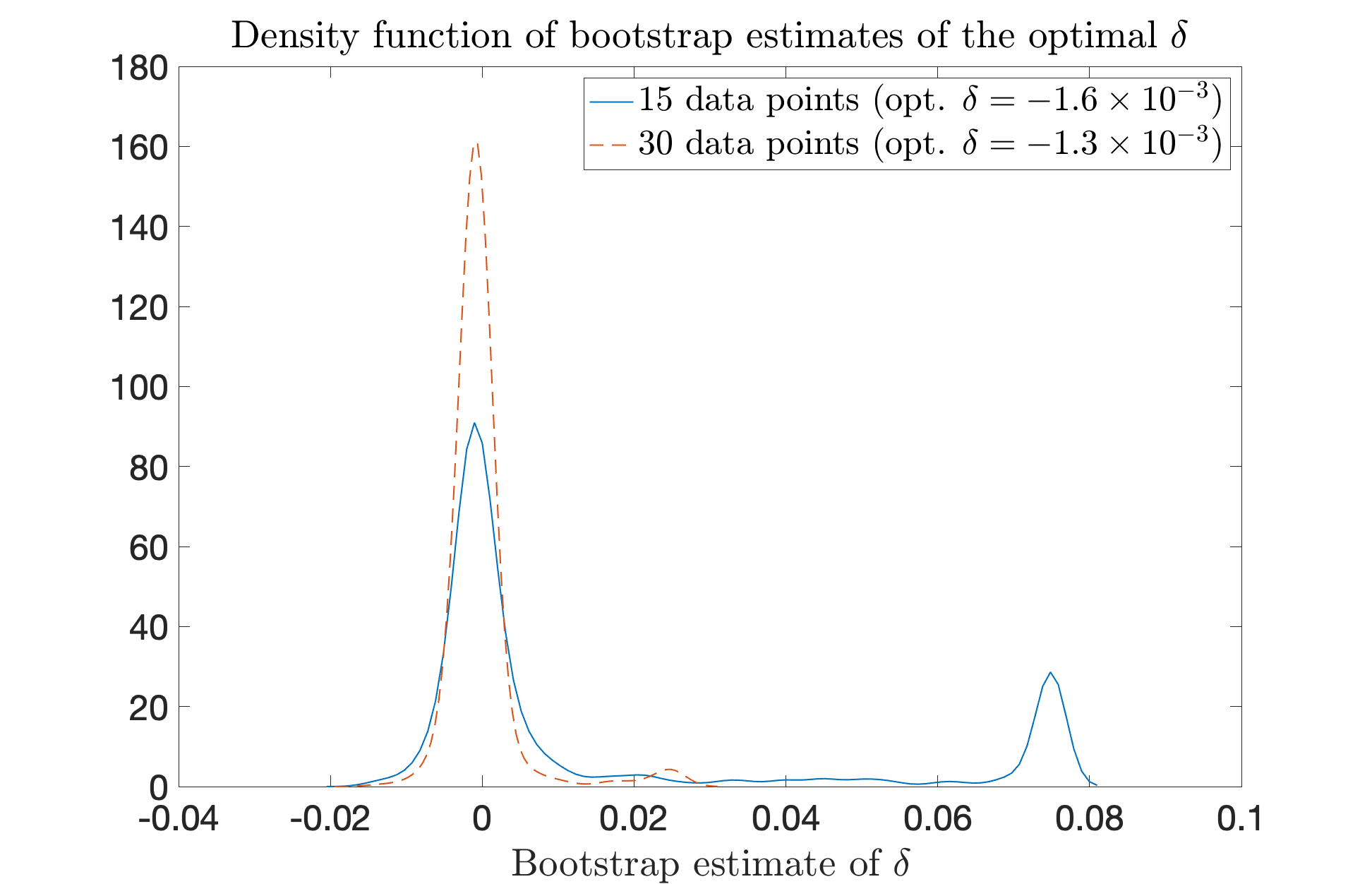}
\caption{Density of bootstrap estimates of $\delta$.}
\label{fig:delta_plots}
\end{subfigure}
\caption{(A) shows the average worst-case sensitivity as a function of $\delta$.
Consistent with \eqref{eq:sensitivity}, it changes linearly around $\delta=0$. (B) shows the (average) in-sample mean-sensitivity frontier corresponding to a range of $\delta$.  The change in the in-sample expected reward around $\delta=0$ (SAA) is small compared to changes in the worst-case sensitivity. (C) shows the out-of-sample mean-sensitivity (variance) frontier that is mapped out when $\delta$ is varied. Out-of-sample expected reward is maximized when $\delta=-1.3\times 10^{-3}$, which is on the optimistic part of the frontier. The expected reward is $193$ which exceeds that of the SAA optimizer (189). (D) shows the distribution of bootstrap estimates of the optimal $\delta$ when there are $n=15$ data points (optimal value $\delta = -1.6\times 10^{-3}$) and $n=30$ data points (optimal value $\delta=-1.3\times 10^{-3}$). Estimation error is large when $n=15$ and there is a long right tail. Accuracy improves significantly when there are $n=30$ data points. }
\label{fig:sensitivity}
\end{figure}



\end{example}

\subsection{Other remarks} \label{sec:other}

The paper \cite{Lam2}\footnote{The paper \cite{Lam2} was posted on 
arXiv one day after the first version of this paper. Both were written independently.} studies data-driven algorithms that admit solutions of the form
\begin{eqnarray}
\tilde{x}_n(\delta) = x^\star(0) + \delta \tilde{\pi} + \frac{1}{\sqrt{n}}\tilde{W}_n + o_P(\delta + n^{-\frac{1}{2}}),
\label{eq:H1}
\end{eqnarray}
where $x^\star(0)$ is the population optimizer \eqref{eq:opt_dgm}, and $\tilde{x}_n(\delta)$ is the data-driven solution when there are $n$ data points, and $\delta$ is the free parameter of the algorithm (e.g. a regularization parameter, the size of an uncertainty set, robustness parameter).

To evaluate performance \eqref{eq:H1}, \cite{Lam2} 
considers the difference (regret) between the expected reward of the data-driven and population optimizers
\begin{eqnarray*}
{\mathcal G}(\tilde{x}_n(\delta)) := g(x^\star(0)) - g(\tilde{x}_n(\delta))
\end{eqnarray*}
where $g(x)$ was defined in \eqref{eq:g}, and  shows that the weak limit ($n\rightarrow\infty$) of the scaled regret $n{\mathcal G}(\tilde{x}_n(\delta_n))$ 
is second-order stochastically larger than the weak limit of the scaled regret $n{\mathcal G}(\tilde{x}_n(0))$  of the SAA optimizer. In this sense it is impossible to beat SAA asymptotically when an algorithm has solutions of the form \eqref{eq:H1}.

While this appears to contradict Theorem \ref{theorem:out-of-sample}, which says that it is always possible to beat SAA using DRO/DOO,
there is actually no inconsistency.
\cite{Lam2} shows the optimality of the {\it limiting reward distribution} of SAA under second order stochastic dominance, whereas we use the asymptotic distribution of the DRO/DOO and SAA solutions to estimate the difference  between the out-of-sample expected rewards for finite $n$ \eqref{eq:out-of-sample-reward}.This difference is of order $\delta/n$, and the advantage of DRO/DOO over SAA vanishes at a rate faster than $1/n$ if $\delta \equiv \delta_n$ itself is going to $0$. Indeed,  it was shown in Theorem \ref{theorem:out-of-sample} that the optimal choice of $\delta$ is of order $1/n$ and the improvement over SAA is of the order $1/n^2$.  If the goal is to beat the expected reward of SAA out-of-sample, the improvement from DRO/DOO is at best modest when data sets are moderate and diminishes quickly in $n$. There is no advantage in the limit, consistent with \cite{Lam2}.

We have argued, however, that the advantage of DRO is more than the possibility that it can sometimes, but not always, beat SAA out-of-sample by just a little bit: It also reduces the sensitivity of the expected reward to misspecification in the nominal which, as the reader may recall, was the point of robust decision making in the first place. More generally, it captures the tradeoff between maximizing expected reward and controlling sensitivity which is relevant in many applications of data-driven decision making.



\section{Conclusion}
Solutions of DRO problems can sometimes have a larger out-of-sample expected reward than SAA. Indeed, much of the literature seems to suggest that this is the aspiration of DRO; uncertainty sets are often calibrated by optimizing an estimate of the out-of-sample expected reward, and recent papers focus on theoretically characterizing the out-of-sample expected reward of DRO solutions; applications of DRO are hailed a success when the out-of-sample mean exceeds that of SAA. However,
beating SAA may not be possible if we only consider worst-case models.

We have shown that it is always possible to beat SAA out-of-sample if we consider both optimistic (DOO) and pessimistic (DRO) perturbations of the nominal model. If the only concern is beating the out-of-sample expected reward of the SAA optimizer, a decision maker should consider best-case and worst-case models.

As tempting as this may sound, there is a catch. While the out-of-sample expected reward might be larger, the worst-case sensitivity of a best-case (DOO) optimizer is larger than that of SAA so it will also be less robust. Indeed, the improvement in the expected reward from using the optimal DOO optimizer is small (of order $n^{-2}$) relative to the robustness cost.
The optimal value of the parameter $\delta$ also needs to be estimated, and bootstrap or cross-validation estimates may be unreliable if the training data set is small.

\clearpage




\newpage

\appendix

\section{First order conditions \eqref{eq:FOC_emp}-\eqref{eq:FOC_pop}}
\label{app-psi}

Consider the population version of the DRO problem \eqref{eq:opt_DRO}
\begin{eqnarray*}
\max_{x, c} G(x, c)
\end{eqnarray*}
where
\begin{eqnarray}
G(\delta, x, c) = -\frac{1}{\delta}{\mathbb E}_{\mathbb P}\Big[\phi^*\Big(-\delta\big[f(x,\,Y)+ c\big]\Big)\Big]-c.
\end{eqnarray}
Differentiating with respect to $x$, the first order conditions are
\begin{eqnarray}
\nabla_x G(\delta, x, c) & = & {\mathbb E}_{\mathbb P}\Big[[\phi^*]'\Big(-\delta\big[f(x, Y) + c\big]\Big)\nabla_xf(x, Y)\Big] = 0, \label{eq:FOC-app}\\ [5pt]
\nabla_c G(\delta, x, c) & = & {\mathbb E}_{\mathbb P}\Big[[\phi^*]'\Big(-\delta\big[f(x, Y) + c\big]\Big)\Big]-1 = 0 \nonumber
\end{eqnarray}
with solutions $(x^\star(\delta), c^\star(\delta))$. We eventually wish to show that $(x^\star(\delta), c^\star(\delta))$ is continuously differentiable in a neighborhood of $\delta=0$, which depends on the properties of the convex conjugate $\phi^*(\zeta)$. Under Assumption \eqref{ass:phi} (see Theorem 3.2 in \cite{GKL}) $\phi^*(\zeta)$ is twice continuously differentiable in the neighborhood of $\zeta=0$ and satisfies
\begin{eqnarray}
\phi^*(\zeta) = \zeta + \frac{1}{2!}\Big(\frac{1}{\phi''(1)}\Big)\zeta^2 +o(\zeta^2).
\label{eq:phi*}
\end{eqnarray}
It follows that $[\phi^*]'(\zeta)$ is continuously differentiable in a neighborhood of $\zeta=0$ and
\begin{eqnarray}
[\phi^*]'(\zeta) = 1 + \frac{\zeta}{\phi''(1)}  + o(\zeta),
\label{eq:phi*-derivative}
\end{eqnarray}
so by \eqref{eq:FOC-app}
\begin{eqnarray*}
\nabla_c G(\delta, x, c) = {\mathbb E}_{\mathbb P}\Big[-\frac{\delta}{\phi''(1)}\big[f(x, Y)+c\big]\Big] + o(\delta)
\end{eqnarray*}
and we can write the first order conditions as
\begin{eqnarray*}
\nabla_x G(\delta, x, c) & = & {\mathbb E}_{\mathbb P}\Big[ [\phi^*]'\Big(-\delta\big[f(x, Y) + c\big]\Big)\nabla_xf(x, Y)\Big] = 0, \\ [5pt]
-\frac{\phi''(1)}{\delta}\nabla_c G(\delta, x, c)  & = &  {\mathbb E}_{\mathbb P}\Big[-\frac{\phi''(1)}{\delta}\Big\{[\phi^*]'\Big(-\delta\big[f(x, Y) + c\big]\Big)-1\Big\}\Big]  = 0.
\end{eqnarray*}
This does not change the solution of \eqref{eq:FOC-app}, but allows us to use the Implicit Function Theorem to establish level of smoothness for $(x^\star(\delta), c^\star(\delta))$ in Proposition \ref{prop:solution_differentiability} (in particular, that \eqref{app-J} is invertible in the proof in Appendix \ref{sec:App-diff}). It follows that $\psi(x, c, Y)$ is given by \eqref{eq:psi}.  First order conditions for the sample DRO problem, and the DOO problems can be derived similarly.

\section{Proof of Proposition \ref{prop:solution_differentiability}} \label{sec:App-diff}
We use the Implicit Function Theorem to show existence of the solution $(x^{\star}(\delta),\,c^{\star}(\delta))$ of \eqref{eq:FOC_pop} and to characterize its smoothness. (Existence and smoothness of $(x_n(\delta),\,c_n(\delta))$ can be shown in a similar way). With a mild abuse of notation, let
\begin{eqnarray*}
g(\delta,\,x, \, c)\equiv \left[\begin{array}{c} g_1(\delta,\,x,\,c) \\ g_2(\delta,\,x,\,c)\end{array}\right] := {\mathbb E}_{\mathbb P}[\psi(x,\,c,\,Y)].
\end{eqnarray*}
The first-order conditions for the population problems \eqref{eq:FOC_pop} are
\begin{eqnarray}
g(\delta,\,x,\,c)= \left[\begin{array}{cc}0 \\ 0\end{array}\right].
\label{foc1}
\end{eqnarray}
Since $[\phi^*]'(\zeta)$ is continuously differentiable in a neighborhood of $\zeta=0$ (Theorem 3.2 in \cite{GKL}), $g(\delta,\,x,\,c)$ is continuously differentiable in $\delta$ in some open neighborhood of $\delta=0$ for every fixed $(x,\,c)$. By \eqref{eq:phi*-derivative}
\begin{eqnarray*}
g(\delta,\,x, \, c)=\left[\begin{array}{c}{\mathbb E}_{{\mathbb P}}\Big[\nabla_x f(x,Y)\Big] -\frac{\delta}{\phi''(1)}{\mathbb E}_{\mathbb P}\Big[ \Big(f(x,Y)+c\Big)\nabla_x f(x,Y)\Big] + o(\delta) \\[10pt] \mathbb{E}_{{\mathbb P}}\Big[f(x,Y)+c \Big]+ O(\delta) \end{array}\right],
\end{eqnarray*}
and
\begin{eqnarray}
g\big(0,\,x^{\star}(0),\,-{\mathbb E}_{{\mathbb P}} [f(x^{\star}(0),\,Y)]\big)=0,
\label{foc1_delta0}
\end{eqnarray}
where $x^{\star}(0)$ is the solution of the empirical problem.
Since $f(x,\,Y)$ is twice continuously differentiable in $x$, $g(\delta,\,x,\,c)$ is continuously differentiable in a neighborhood of $\big(0,\,x^{\star}(0),\,-{\mathbb E}_{{\mathbb P}} [f(x^{\star}(0),\,Y)]\big)$, and
\begin{eqnarray}
J_{g,\,(x,\,c)}(\delta,\,x,\,c)\Big|_{(0,\,x^{\star}(0),\,-{\mathbb E}_{{\mathbb P}} [f(x^{\star}(0),\,Y)])}
 & \equiv & \left.\left[\begin{array}{cc}\nabla_x g_1(\delta,\,x,\,c) & \nabla_c g_1(\delta,\,x,\,c) \nonumber \\
\nabla_x g_2(\delta,\,x,\,c) & \nabla_c g_2(\delta,\,x,\,c)
\end{array}\right]\right|_{(0,\,x^{\star}(0),\, -{\mathbb E}_{{\mathbb P}} [f(x^{\star}(0),\,Y)])}\\ & = & \left[\begin{array}{cc}{\mathbb E}_{\mathbb P}[\nabla_x^2f(x^{\star}(0),\,Y)] & 0 \\ 0 & 1\end{array}\right]\label{app-J}
\end{eqnarray}
is invertible. It follows from the Implicit Function Theorem that $(x^{\star}(\delta),\,c^{\star}(\delta))$ exists and is continuously differentiable in an open neighborhood of $\big(0,\,x^{\star}(0),\,-{\mathbb E}_{{\mathbb P}} [f(x^{\star}(0),\,Y)]\big)$ so we can write
\begin{eqnarray*}
x^\star(\delta) & = & x^\star(0) + \delta \pi + o(\delta) \\ [1pt]
c^\star(\delta) & = & - {\mathbb E}_{\mathbb P}[f(x^\star(0), Y)] + \delta \kappa + o(\delta)
\end{eqnarray*}
where $\pi \in {\mathbb R}^d$ and $\kappa$ is a constant. Substituting into the first equation in \eqref{foc1_delta0}, a Taylor series expansion around $\delta=0$ gives
\begin{eqnarray*}
\lefteqn{{\mathbb E}_{\mathbb P}\big[\nabla_x f(x^\star(0), Y)\big] + \delta \pi' {\mathbb E}_{\mathbb P}\big[\nabla^2_x f(x^\star(0), Y)\big]} \\ [5pt]
& &  - \frac{\delta}{\phi^{''}(1)}{\mathbb E}_{\mathbb P}\Big[\Big(f(x^\star(0), Y) - {\mathbb E}_{\mathbb P}[f(x^\star(0), Y)\Big) \nabla_x f(x^\star(0), Y) \Big] + o(\delta) = 0.
\end{eqnarray*}
Optimality of $x^\star(0)$ for the population problem implies ${\mathbb E}_{\mathbb P}\big[\nabla_x f(x^\star(0), Y)\big] =0$ while the coefficient of the $\delta$ term is zero if
\begin{eqnarray*}
\pi =  \frac{1}{\phi^{''}(1)} \Big({\mathbb E}_{\mathbb P}\big[\nabla^2_x f(x^\star(0), Y)\big]\Big)^{-1} \mbox{Cov}_{\mathbb P}\Big(f(x^\star(0), Y), \nabla_x f(x^\star(0), Y)\Big).
\end{eqnarray*}
Smoothness of $(x_n(\delta),\,c_n(\delta))$ and the expression \eqref{eq:insample_rob_asymp_bias}-\eqref{eq:insample_pi} can be established in a similar manner.






\section{Proof of Proposition \ref{prop:distribution-of-solution}} \label{app:prop-stat}

We begin with some preliminaries. Let $\psi(x, c, Y)$ be defined by \eqref{eq:psi}, and
\begin{eqnarray*}
\nabla\psi(x^{\star}(\delta),\,c^{\star}(\delta), Y)
= \left[
\begin{array}{ccc|c}
\nabla_{x_1} \psi_1 & \ldots & \nabla_{x_m} \psi_1 & \nabla_c \psi_1 \\
\hline
\nabla_{x_1} \psi_2 & \ldots & \nabla_{x_m} \psi_2 & \nabla_c \psi_2 \\
\end{array}
\right](x^{\star}(\delta),\,c^{\star}(\delta), Y)
\end{eqnarray*}
be the Jacobian of $\psi(x,c, Y)$ 
and the matrices
\begin{eqnarray*}
\widetilde{A}(\delta) & = &  {\mathbb E}_{\mathbb{P}}[\nabla\psi(x^{\star}(\delta),\,c^{\star}(\delta), Y)] \in {\mathbb R}^{(d+1)\times (d+1)},\\
\widetilde{B}(\delta)  & = &  {\mathbb E}_{\mathbb P}[\psi(x^\star(\delta),\,c^\star(\delta), Y)\,\psi(x^\star(\delta),\,c^\star(\delta), Y)'] \in {\mathbb R}^{(d+1) \times (d+1)}.
\end{eqnarray*}
We also define the ${\mathbb R}^{d+1}$ valued random vectors
\begin{eqnarray}
\widetilde{W}_n(\delta) & = & -\widetilde{A}(\delta)^{-1}\frac{1}{\sqrt{n}}\sum_{i=1}^n\psi(x^\star(\delta),\,c^\star(\delta), \, Y_i), \nonumber \\
\widetilde{V}_n(\delta) & = & \widetilde{H}_n(\delta) \widetilde{W}_n(\delta) - \widetilde{I}_n(\delta)
\label{eq:VW}
\end{eqnarray}
where
\begin{eqnarray*}
\widetilde{H}_n(\delta) & = & -\frac{1}{\sqrt{n}}\widetilde{A}(\delta)^{-1}\sum_{i=1}^n\Big\{ \nabla\psi(x^{\star}(\delta),\,c^{\star}(\delta), Y_i) - {\mathbb E}_{\mathbb P}\big[\nabla\psi(x^{\star}(\delta),\,c^{\star}(\delta), Y)\big]\Big\} \\
\tilde{I}_n(\delta) & = & \widetilde{A}(\delta)^{-1}
\left[\begin{array}{c}
\widetilde{W}_n(\delta)'{\mathbb E}_{\mathbb P}\big[\nabla^2\psi^{(1)}(x^\star(\delta),\,c^\star(\delta), \, Y) \big]\widetilde{W}_n(\delta)  \\
\vdots \\
\widetilde{W}_n(\delta)'{\mathbb E}_{\mathbb P}\big[\nabla^2\psi^{(d+1)}(x^\star(\delta),\,c^\star(\delta), \, Y)\big]\widetilde{W}_n(\delta)
\end{array}\right]
\end{eqnarray*}
Here, $\psi^{(i)}(x, c, Y)$ is the $i^{th}$ component of the $d+1$ dimensional function
\begin{eqnarray*}
\psi(x, c, Y) = \left[\begin{array}{c}\psi^{(1)}(x, c, Y) \\ \vdots \\ \psi^{(d+1)}(x, c, Y)\end{array}\right]
\end{eqnarray*}
that was defined in \eqref{eq:psi}, and
$\nabla^2\psi^{(i)}(x^\star(\delta),\,c^\star(\delta), \, Y)$
is the Hessian $\nabla^2 \psi^{(i)}(x,\,c, \, Y)$ evaluated at $(x^\star(\delta), c^\star(\delta))$. Note that condition (2) of Assumption \ref{ass:reg} implies that $\phi^*(\zeta)$ is three times continuously differentiable\footnote{It can be shown, along the lines of the Proof of Theorem 3.2 in \cite{GKL}, that
\begin{eqnarray*}
\phi^*(\zeta) = \zeta + \frac{\zeta^2}{2!}\Big(\frac{1}{\phi^{(2)}(1)}\Big) - \frac{\zeta^3}{3!}\Big(\frac{1}{\phi^{(3)}(1)}\frac{1}{\phi^{(2)}(1)}\Big) + o(\zeta^3)
\end{eqnarray*}
where $\phi^{(i)}(z)$ is the $i^{th}$ derivative of $\phi(z)$ with respect to $z$. } in a neighborhood of $\zeta=0$. Together with condition (1), it follows that  the Hessian $\psi^{(i)}(x,\,c, \, Y)$ evaluated at $(x^\star(\delta), c^\star(\delta))$ above is well defined.

The following result characterizes the distribution of $(x_n(\delta), \, c_n(\delta))$ and is obtained by applying  Theorem 5.21 in \cite{vdV} and Lemma 1 from \cite{KR} to \eqref{eq:psi}.
\begin{proposition} \label{prop:vdV}
Suppose that data $\{Y_1,\cdots,\,Y_n\}$ is drawn iid from $\mathbb P$, $f(x,\,Y)$ satisfies Assumption \ref{ass:f},  $\phi(z)$ satisfies Assumption \ref{ass:phi}, and 
 $\delta$ is such that Assumption \ref{ass:reg} holds.
Then there is a unique solution $(x^\star(\delta),\,c^\star(\delta))$ of the first-order conditions \eqref{eq:FOC_pop}, the matrix $\widetilde{A}(\delta)$ is invertible, $(x_n(\delta),\,c_n(\delta)) \overset{P}{\longrightarrow} (x^\star(\delta),\,c^\star(\delta))$, and for $n$ sufficiently large
\begin{eqnarray*}
\left[\begin{array}{c}x_n(\delta)\\ c_n(\delta) \end{array}\right] = \left[\begin{array}{c}x^\star(\delta) \\ c^\star(\delta)\end{array}\right] + \frac{1}{\sqrt{n}} \widetilde{W}_n(\delta) + \frac{1}{n}\widetilde{V}_n(\delta)+ o_P(n^{-3/2}).
\label{eq:xc}
\end{eqnarray*}
$\widetilde{W}_n(\delta)$ has mean $0$ and covariance matrix
\begin{eqnarray*}\widetilde{\xi}(\delta) \equiv {\mathbb V}_{\mathbb P}[\widetilde{W}_n(\delta)]= \widetilde{A}(\delta)^{-1} \widetilde{B}(\delta) {\widetilde{A}(\delta)^{-1}}',\end{eqnarray*}
and $\widetilde{V}_n(\delta)$ is $O_P(1)$ with mean
\begin{eqnarray*}
{\mathbb E}_{\mathbb P}[\widetilde{V}_n(\delta)]
& = & {\mathbb E}_{\mathbb P} \Big[\widetilde{A}(\delta)^{-1} \nabla\psi(x^{\star}(\delta),\,c^{\star}(\delta), Y)\widetilde{A}(\delta)^{-1}\psi(x^{\star}(\delta),\,c^{\star}(\delta), Y)\Big] \\ & &
- \widetilde{A}(\delta)^{-1}
\left[\begin{array}{c}
{\rm tr} \big\{\widetilde{\xi}(\delta) \, {\mathbb E}_{\mathbb P}\big[\nabla^2\psi^{(1)}(x^\star(\delta),\,c^\star(\delta), \, Y) \big] \big\} \\
\vdots \\
{\rm tr} \big\{\widetilde{\xi}(\delta)\,{\mathbb E}_{\mathbb P}\big[\nabla^2\psi^{(d+1)}(x^\star(\delta),\,c^\star(\delta), \, Y)\big]\big\}
\end{array}\right].
\end{eqnarray*}
\end{proposition}

\medskip

Proposition \ref{prop:distribution-of-solution} is largely a restatement of Proposition \ref{prop:vdV} for the $x_n(\delta)$ component of $(x_n(\delta), c_n(\delta))$. In particular,  expression \eqref{eq:x1} for $x_n(\delta)$ holds because we can extract the ${\mathbb R}^d$ valued random vectors $W_n(\delta)$ and $V_n(\delta)$ of \eqref{eq:VW}
 corresponding to $x_n(\delta)$ in \eqref{eq:xc}
\begin{eqnarray*}
\widetilde{W}_n(\delta)  =  \left[\begin{array}{c}W_n(\delta) \\ U_n(\delta)\end{array}\right] \in {\mathbb R}^{d +1}, \hspace{0.5cm}
\widetilde{V}_n(\delta)  =  \left[\begin{array}{c}V_n(\delta) \\ T_n(\delta)\end{array}\right]\in {\mathbb R}^{d+1},
\end{eqnarray*}
and  $\xi(\delta)\in {\mathbb R}^{d\times d}$ and $\overline{V}(\delta)\in{\mathbb R}^d$ from
\begin{eqnarray*}
\widetilde{\xi}(\delta) = \left[\begin{array}{cc} \xi(\delta) &
\kappa(\delta) \\
\kappa(\delta)' & \eta(\delta)\end{array}\right], \hspace{0.5cm} {\mathbb E}_{\mathbb P}[\widetilde{V}_n(\delta)]= \left[\begin{array}{c}\overline{V}(\delta) \\ \overline{S}(\delta)\end{array}\right].
\end{eqnarray*}
$\xi(\delta)$ and $\overline{V}(\delta)$ are continuously differentiable in a neighborhood of $\delta=0$ because the pair $(x^\star(\delta), c^\star(\delta))$ is continuously differentiable in a neighborhood of $\delta=0$ (Proposition \ref{prop:solution_differentiability}) and the conditions imposed on $\psi(x, c, Y)$ in Assumption \ref{ass:reg}.

\section{Assumptions \ref{ass:SAALp} and \ref{ass:DLp}} \label{app:Lp}


We know from \eqref{eq:x0} that the solution of the SAA problem has an expansion of the form
\begin{eqnarray}
x_n(0) = x^\star(0) + \frac{1}{\sqrt{n}} W_n(0)  + \frac{1}{n} V_n(0) + H_n(0)
\label{sol:exp}
\end{eqnarray}
where the error term $H_n(0)$ is $O_p(n^{-\frac{3}{2}})$. Assumption \ref{ass:SAALp} holds if $\sup_n {\mathbb E}_{\mathbb P}\|n^\frac{3}{2}{H}_n(0)\|^p<\infty$ for $p=1, 2$. We now derive an expression for $H_n(0)$ to determine conditions under which this holds. 
To ease notation, we assume $x$ is scalar. Our derivation shows that $n^{\frac{3}{2}} H_n(\delta)$ is a product of derivatives of $f(x, Y)$ with respect to $x$ of up to order $4$, and Assumption \ref{ass:SAALp} holds if the first two moments of $n^{\frac{3}{2}} H_n(\delta)$ are finite. 

Our approach  generalizes to the case where the decision variable is multi-dimensional, which is required for DRO and DOO (Assumption \ref{ass:DLp}) where the pair $(x, c)$ is the decision and SAA with vector-valued $x$, at the cost of even uglier equations and tensor notation. Not surprisingly, the same insight, that products of (partial) derivatives of order up to 4 have finite first- and second moment, holds. The equations are long but not particularly insightful (beyond this observation) and the derivation is arduous. We outline how this can be done and leave details to the reader.

To derive the expansion \eqref{sol:exp}, we embed the first order conditions for our sample problem in a family of fixed-point equations $G(y, \varepsilon)=0$ parameterized by $\varepsilon\in[0, 1]$ where
\begin{eqnarray*}
G(y, \varepsilon) & = & {\mathbb E}_{\mathbb P}[\nabla_xf(y, Y)] + \varepsilon\frac{1}{n}\sum_{i=1}^n\Big\{\nabla_xf(y, Y_i) - {\mathbb E}_{\mathbb P}[\nabla_xf(y, Y)]\Big\}\\
& = & {\mathbb E}_{\mathbb P}[\nabla_xf(y, Y)] + \varepsilon\frac{1}{n}\sum_{i=1}^n \nabla_x \widetilde{f}(y, Y_i).
\end{eqnarray*}
Here
\begin{eqnarray*}
\widetilde{f}(y, Y) = f(y, Y) - {\mathbb E}_{\mathbb P}[f(y, Y)]
\end{eqnarray*}
is the deviation of the random reward $f(y, Y)$ from its expectation ${\mathbb E}_{\mathbb P}[f(y, Y)]$, while $\nabla_x \widetilde{f}(y, Y)$ is the derivative of the centered random variable. We denote the solution of $G(y, \varepsilon)=0$ as $y(\varepsilon)$. 
(We change the argument in $G$ from $x$ to $y$ to minimize the chance of confusion between $x_n(0)$, the solution of the SAA problem, and the solution $y(\varepsilon)$ of the fixed point problem).
The case $\varepsilon=0$ 
\begin{eqnarray*}
G(y(0), 0) = {\mathbb E}_{\mathbb P}[\nabla_xf(y(0), Y)] =0
\end{eqnarray*}
corresponds to the population problem with solution $y(0) = x^\star(0)$; $\varepsilon=1$ is the SAA problem with solution $y(1) = x_n(0)$
\begin{eqnarray*}
G(x_n(0),1) = \frac{1}{n}\sum_{i=1}^n \nabla_xf(x_n(0), Y_i) = 0.
\end{eqnarray*}

To begin, we find the derivatives of $y(\varepsilon)$ with respect to $\varepsilon$, or $A(\varepsilon)$, $B(\varepsilon)$ and $C(\varepsilon)$, such that
\begin{eqnarray*}
y(\varepsilon+\Delta) = y(\varepsilon) + \Delta A(\varepsilon) + \frac{\Delta^2}{2}B(\varepsilon) + \frac{\Delta^3}{3!}C(\varepsilon) + o(\Delta^3).
\end{eqnarray*}
Since $G(y(\varepsilon+\Delta), \varepsilon + \Delta)=0$ we do this by expanding $G(y(\varepsilon+\Delta), \varepsilon + \Delta)$ around $G(y(\varepsilon), \varepsilon) $ and choosing $A(\varepsilon)$, $B(\varepsilon)$ and $C(\varepsilon)$  to ensure the first-order condition holds. 

Specifically, a 
Taylor series expansion around $\varepsilon$ gives (after some work)
\begin{eqnarray*}
\lefteqn{G(y(\varepsilon+\Delta), \varepsilon + \Delta) = G(y(\varepsilon), \varepsilon)  }\\ [5pt]
& & + \Delta\left\{\Big[{\mathbb E}_{\mathbb P}[\nabla_x^2f(y(\varepsilon), Y)] + \varepsilon\frac{1}{n}\sum_{i=1}^n \nabla_x^2 \widetilde{f}(y(\varepsilon), Y_i)\Big] A(\varepsilon) + \frac{1}{n}\sum_{i=1}^n \nabla_xf(y(\varepsilon), Y_i)\right\} \\ [5pt]
& & + \frac{\Delta^2}{2!}\left\{B(\varepsilon) \Big({\mathbb E}_{\mathbb P}[\nabla_x^2f(y(\varepsilon), Y)] + \varepsilon\frac{1}{n}\sum_{i=1}^n \nabla_x^2 \widetilde{f}(y(\varepsilon), Y_i) \Big)\right. \\
& & \left.+ A(\varepsilon)^2 \Big({\mathbb E}_{\mathbb P}[\nabla_x^3f(y(\varepsilon), Y)] + \varepsilon\frac{1}{n}\sum_{i=1}^n \nabla_x^3 \widetilde{f}(y(\varepsilon), Y_i)\Big)+  A(\varepsilon) \frac{1}{n}\sum_{i=1}^n2 \nabla^2_x \widetilde{f}(y(\varepsilon), Y)\right\} \\
& & + \frac{\Delta^3}{3!}\left\{C(\varepsilon)\Big({\mathbb E}_{\mathbb P}[\nabla_x^2f(y(\varepsilon), Y)] + \varepsilon\frac{1}{n}\sum_{i=1}^n \nabla_x^2 \widetilde{f}(y(\varepsilon), Y_i) \Big)  \right.\\ [5pt]
& &  + A(\varepsilon)^3 \Big({\mathbb E}_{\mathbb P}[\nabla_x^4f(y(\varepsilon), Y)] + \varepsilon\frac{1}{n}\sum_{i=1}^n \nabla_x^4 \widetilde{f}(y(\varepsilon), Y_i) \Big)  \\
& &  + 3 A(\varepsilon) B(\varepsilon) \Big({\mathbb E}_{\mathbb P}[\nabla_x^3f(y(\varepsilon), Y)] + \varepsilon\frac{1}{n}\sum_{i=1}^n \nabla_x^3 \widetilde{f}(y(\varepsilon), Y_i) \Big) \\
& & \left.+ 3 \Big(\frac{1}{n}\sum_{i=1}^n\nabla^2_x\widetilde{f}(y(\varepsilon), Y_i) \Big) B(\varepsilon) +  3 \Big(\frac{1}{n}\sum_{i=1}^n\nabla^3_x\widetilde{f}(y(\varepsilon), Y_i) \Big) A(\varepsilon)^2\right\} + o(\Delta^3).
\end{eqnarray*} 
Since $G(y(\varepsilon+\Delta), \varepsilon+\Delta)=0$,  the coefficients of the powers of $\Delta$ must be zero and hence
\begin{eqnarray*}
A(\varepsilon)   =   \frac{1}{\sqrt{n}} \Lambda_n(\varepsilon), \hspace{0.5cm} B(\varepsilon)   =   \frac{1}{{n}} \Gamma_n(\varepsilon), \hspace{0.5cm}
C(\varepsilon)  =  n^{-\frac{3}{2}} \Theta_n(\varepsilon)
\end{eqnarray*}
where
\begin{eqnarray*}
\Lambda_n(\varepsilon) & = & \frac{-\frac{1}{\sqrt{n}}\sum_{i=1}^n \nabla_xf(y(\varepsilon), Y_i)}{{\mathbb E}_{\mathbb P}[\nabla_x^2f(y(\varepsilon), Y)] + \varepsilon\frac{1}{n}\sum_{i=1}^n \nabla_x^2 \widetilde{f}(y(\varepsilon), Y_i)}\\ [5pt]
\Gamma_n(\varepsilon) & = &- [\Lambda_n(\varepsilon)]^2 \; \left\{ \frac{ {\mathbb E}_{\mathbb P}[\nabla_x^3f(y(\varepsilon), Y)] + \varepsilon\frac{1}{n}\sum_{i=1}^n \nabla_x^3 \widetilde{f}(y(\varepsilon), Y_i)}{{\mathbb E}_{\mathbb P}[\nabla_x^2f(y(\varepsilon), Y)] + \varepsilon\frac{1}{n}\sum_{i=1}^n \nabla_x^2 \widetilde{f}(y(\varepsilon), Y_i)} \right\} \\ [5pt]
 & &  -   \Lambda_n(\varepsilon)\left\{\frac{\frac{1}{\sqrt{n}}\sum_{i=1}^n2 \nabla^2_x \widetilde{f}(y(\varepsilon), Y)}{{\mathbb E}_{\mathbb P}[\nabla_x^2f(y(\varepsilon), Y)] + \varepsilon\frac{1}{n}\sum_{i=1}^n \nabla_x^2 \widetilde{f}(y(\varepsilon), Y_i)}\right\} \\
 \Theta_n(\varepsilon) & = & -3 {\Gamma_n(\varepsilon)}\left\{\frac{\frac{1}{\sqrt n}\sum_{i=1}^n\nabla^2_x\widetilde{f}(y(\varepsilon), Y_i)}{{\mathbb E}_{\mathbb P}[\nabla_x^2f(y(\varepsilon), Y)] + \varepsilon\frac{1}{n}\sum_{i=1}^n \nabla_x^2 \widetilde{f}(y(\varepsilon), Y_i)}\right\} \\ [5pt]
& & -3 [\Lambda_n(\varepsilon)]^2\left\{\frac{\frac{1}{\sqrt n}\sum_{i=1}^n\nabla^3_x\widetilde{f}(y(\varepsilon), Y_i)}{{\mathbb E}_{\mathbb P}[\nabla_x^2f(y(\varepsilon), Y)] + \varepsilon\frac{1}{n}\sum_{i=1}^n \nabla_x^2 \widetilde{f}(y(\varepsilon), Y_i)}\right\} \\ [5pt]
& & - [\Lambda_n(\varepsilon)]^3\left\{\frac{{\mathbb E}_{\mathbb P}[\nabla_x^4f(y(\varepsilon), Y)] + \varepsilon\frac{1}{n}\sum_{i=1}^n \nabla_x^4 \widetilde{f}(y(\varepsilon), Y_i)}{{\mathbb E}_{\mathbb P}[\nabla_x^2f(y(\varepsilon), Y)] + \varepsilon\frac{1}{n}\sum_{i=1}^n \nabla_x^2 \widetilde{f}(y(\varepsilon), Y_i)} \right\} \\ [5pt]
& & - 3 \Lambda_n(\varepsilon) \Gamma_n(\varepsilon) \left\{\frac{{\mathbb E}_{\mathbb P}[\nabla_x^3f(y(\varepsilon), Y)] + \varepsilon\frac{1}{n}\sum_{i=1}^n \nabla_x^3 \widetilde{f}(y(\varepsilon), Y_i)}{{\mathbb E}_{\mathbb P}[\nabla_x^2f(y(\varepsilon), Y)] + \varepsilon\frac{1}{n}\sum_{i=1}^n \nabla_x^2 \widetilde{f}(y(\varepsilon), Y_i)}\right\}
\end{eqnarray*}
which gives us the derivatives of $y(\varepsilon)$ we are looking for.
This gives the Taylor series expansion
\begin{eqnarray*}
y(\varepsilon) = y(0) + \varepsilon A(0) +\frac{\varepsilon^2}{2} B(0)  + \frac{\varepsilon^{3}}{3!} C(\eta)
\end{eqnarray*}
where $\eta$ is a random variable taking values in $[0, 1]$ that depends on the data. In particular, with $\varepsilon=1$ and noting that $y(\varepsilon) = x_n(0)$ when $\varepsilon=1$
\begin{eqnarray*}
x_n(0) = x^\star(0) + \frac{1}{\sqrt{n}} \Lambda_n(0)  + \frac{1}{n} \Gamma_n(0) + n^{-\frac{3}{2}}\Theta_n(\eta)
\end{eqnarray*}
 where
\begin{eqnarray*}
\Lambda_n(0) & =& \frac{-\frac{1}{\sqrt{n}}\sum_{i=1}^n \nabla_xf(x^\star(0), Y_i)}{{\mathbb E}_{\mathbb P}[\nabla_x^2f(x^\star(0), Y)]} \\
\Gamma_n(0) & = &  - [\Lambda_n(0)]^2 \left\{\frac{ {\mathbb E}_{\mathbb P}[\nabla_x^3f(x^\star(0), Y)] }{{\mathbb E}_{\mathbb P}[\nabla_x^2f(x^\star(0), Y)]} \right\}   -   \Lambda_n(0)\left\{\frac{\frac{1}{\sqrt{n}}\sum_{i=1}^n2 \nabla^2_x \widetilde{f}(x^\star(0), Y)}{{\mathbb E}_{\mathbb P}[\nabla_x^2f(x^\star(0), Y)]}\right\},
\end{eqnarray*}
and $\eta$ is a random variable taking values between $0$ and $1$. Observe that $\Lambda_n(0) = W_n(0)$ and $\Gamma_n(0) = V_n(0)$, as defined in \eqref{eq:x0} and that the error term in \eqref{sol:exp} is $H_n(0) = n^{-\frac{3}{2}}\Theta_n(\eta)$. Assumption \ref{ass:SAALp}  holds if ${\mathbb E}_{\mathbb P}\|\Theta_n(\eta)\|^p<\infty$ for $p=1, 2$. This will be the case if there is a constant $\gamma>0$ such that  $\nabla^2_x f(x, Y) > \gamma$ for every $(x, Y)$ and the derivatives $\|\nabla_x^k f(x, Y)\|$ of order $k =1, 2, 3, 4$ are uniformly bounded, though this is not the mildest assumption.

In the case of DRO/DOO, we have first-order conditions
\begin{eqnarray}
 {\mathbb E}_{\mathbb P} \left[\begin{array}{c}
  [\phi^*]'\Big(-\delta\big[f(x, Y) + c\big]\Big)\nabla_xf(x, Y) = 0 \\ [10pt]
 -\frac{\phi''(1)}{\delta}\Big\{[\phi^*]'\Big(-\delta\big[f(x, Y) + c\big]\Big)-1\Big\}
 \end{array}\right] = 0
 \label{DRO-pop}
\end{eqnarray}
for the solution $(x^\star(\delta), c^\star(\delta))$ of the population problem, and 
\begin{eqnarray}
\frac{1}{n} \sum_{i=1}^n   \left[\begin{array}{c}
  [\phi^*]'\Big(-\delta\big[f(x, Y) + c\big]\Big)\nabla_xf(x, Y) \\ [10pt]
 -\frac{\phi''(1)}{\delta}\Big\{[\phi^*]'\Big(-\delta\big[f(x, Y) + c\big]\Big)-1\Big\}
 \end{array}\right] = 0
 \label{DRO-emp}
\end{eqnarray}
for the solution $(x^\star_n(\delta), c^\star_n(\delta))$ of the empirical version. Analogous to the case of SAA, we can define $(y(\varepsilon), \omega(\varepsilon))$ as the solution of $G((y(\varepsilon), \omega(\varepsilon)), \varepsilon)=[0, 0]'$ ($\varepsilon\in [0, 1]$) where
\begin{eqnarray*}
G((y(\varepsilon), \omega(\varepsilon)), \varepsilon)&  = & (1-\varepsilon) {\mathbb E}_{\mathbb P} \left[\begin{array}{c}
  [\phi^*]'\Big(-\delta\big[f(y(\varepsilon), Y) + \omega(\varepsilon)\big]\Big)\nabla_xf(y(\varepsilon), Y)  \\ [10pt]
 -\frac{\phi''(1)}{\delta}\Big\{[\phi^*]'\Big(-\delta\big[f(y(\varepsilon), Y) + c(\varepsilon)\big]\Big)-1\Big\}
 \end{array}\right]  \\ & &  +  \varepsilon \frac{1}{n} \sum_{i=1}^n \left[\begin{array}{c}
  [\phi^*]'\Big(-\delta\big[f(y(\varepsilon), Y) + \omega\big]\Big)\nabla_xf(y(\varepsilon), Y)
  \\ [10pt]
 -\frac{\phi''(1)}{\delta}\Big\{[\phi^*]'\Big(-\delta\big[f(y(\varepsilon), Y) + \omega(\varepsilon)\big]\Big)\Big\} 
 \end{array}\right].
\end{eqnarray*}
This allows us to embed \eqref{DRO-pop} and \eqref{DRO-emp} in the problem of solving $G((y(\varepsilon), \omega(\varepsilon)), \varepsilon)=[0, 0]'$. Expressions for $\Lambda_n(\varepsilon)$, $\Gamma_n(\varepsilon)$ and the remainder term $\Theta_n(\varepsilon)$ such that
\begin{eqnarray*}
\left[\begin{array}{c}y(\varepsilon)\\\omega(\varepsilon)\end{array}\right] = \frac{1}{\sqrt{n}}\Lambda_n(\varepsilon) + \frac{1}{n} \Gamma_n(\varepsilon) + n^{-\frac{3}{2}}\Theta_n(\eta),
\end{eqnarray*}
where $\eta$ is a random variable between $0$ and $\varepsilon$, can be derived, though equations are long and unwieldy. Setting $\varepsilon=1$ we get
\begin{eqnarray*}
\left[\begin{array}{c}x_n(\delta)\\ c_n(\delta)\end{array}\right] = \frac{1}{\sqrt{n}}W_n(\delta) + \frac{1}{n} V_n(\delta) + n^{-\frac{3}{2}}\Theta_n(\eta),
\end{eqnarray*}
where $W_n(\delta) = \Lambda_n(1)$, $V_n(\delta) = \Gamma_n(1)$ and the
remainder term $\Theta_n(\eta)$ is a product of partial derivatives of 
\begin{eqnarray*} 
\Phi(y, \omega, Y) := -\frac{1}{\delta}\phi^* \Big(-\delta \Big[f(y, Y) + \omega\Big]\Big) - \omega
\end{eqnarray*} 
with respect to $y$ and $\omega$ of up to order $4$, evaluated at $y(\eta)$ for some random $\eta\in [0, 1]$ that depends on the data. Assumption \ref{ass:DLp} holds if $\sup_n {\mathbb E}_{\mathbb P}\|\Theta_n(\eta)\|^p < \infty$ for $p=1, 2$. Again, this will be the case if, for example, the partial derivatives up to order $4$ are uniformly bounded.

As a final note, \cite{Nishiyama} provides conditions for the error term in an $M$-estimation problem  to converge in $L^p$ for every $p\geq 1$. It can also be used to derive conditions for $\sup_n{\mathbb E}_{\mathbb P}\|n^{\frac{3}{2}}H_n\|^p<\infty$  for every $p\geq 1$, and hence for Assumptions \ref{ass:SAALp} and \ref{ass:DLp} to hold.

\section{Proof of Proposition \ref{prop:out-of-sample reward}}
\label{App:DRO-obj}

Recall from  \eqref{eq:bias-var-SAA} that
\begin{eqnarray}
{\mathbb E}_{\mathbb P}[x_n(0)]  = x^\star(0) + \frac{1}{n} \overline{V}(0) + o(n^{-1})
\label{eq:bias-var-SAA-app}
\end{eqnarray}
and hence
\begin{eqnarray}
\nabla_x g({\mathbb E}_{\mathbb P}[x_n(0)]) & = & \nabla_x g\Big(x^\star(0) + \frac{1}{n}\overline{V}(0)+o(n^{-1})\Big) \nonumber \\
& = & \nabla_x g(x^\star(0)) + \nabla^2_x g(x^\star(0))\Big(\frac{1}{n}\overline{V}(0)\Big) + o(n^{-1}) \nonumber \\
& = &  \frac{1}{n} \nabla^2_x g(x^\star(0))\overline{V}(0) + o(n^{-1}) \nonumber \\ [8pt]
\nabla^2_x g({\mathbb E}_{\mathbb P}[x_n(0)]) & = & \nabla^2_x g(x^\star(0)) + O(n^{-1}).
\label{eq:expg_app}
\end{eqnarray}
From Proposition \ref{prop:distribution-of-solution}, the DRO/DOO solution
\begin{eqnarray*}
x_n(\delta) = x^\star(\delta) + \frac{1}{\sqrt{n}} W_n(\delta) + \frac{1}{n}V_n(\delta) + o_P(n^{-1})
\end{eqnarray*}
where $x^\star(\delta) = x^\star(0) + \delta \pi + o(\delta)$. Taking expectations gives
\begin{eqnarray}
 {\mathbb E}_{\mathbb P}[x_n(\delta)] & = & x^\star(\delta)+\frac{1}{n}\overline{V}(\delta)  + o(n^{-1}) \label{eq:Ex-delta-app}\\
 {\mathbb V}_{\mathbb P}[x_n(\delta)] & = & \frac{1}{n}\xi(\delta) + o(n^{-1}).\nonumber
\end{eqnarray}
Using \eqref{eq:bias-var-SAA-app} we can write \eqref{eq:Ex-delta-app} as
\begin{eqnarray}
{\mathbb E}_{\mathbb P}[x_n(\delta)] & = & {\mathbb E}_{\mathbb P}[x_n(0)] + [x^\star(\delta)-x^\star(0)] + \frac{1}{n}[\overline{V}(\delta)-\overline{V}(0)]+ o(n^{-1}) \nonumber \\
&= & {\mathbb E}_{\mathbb P}[x_n(0)] + \delta \pi + o(\delta) + \frac{1}{n}\{\delta \overline{V}_\delta(0) + o(\delta)\} + o(n^{-1}).
\label{eq:exp-Ex-app}
\end{eqnarray}

We now expand $g({\mathbb E}_{\mathbb P}[x_n(\delta)])$ around $g({\mathbb E}_{\mathbb P}[x_n(0)])$:
\begin{eqnarray*}
\lefteqn{g({\mathbb E}_{\mathbb P}[x_n(\delta)])} \\
& = & g({\mathbb E}_{\mathbb P}[x_n(0)]) + \Big(\nabla_x g\big({\mathbb E}_{\mathbb P}[x_n(0)]\big)\Big)'\Big\{{\mathbb E}_{\mathbb P}[x_n(\delta)]-{\mathbb E}_{\mathbb P}[x_n(0)]\Big\} \\
& & + \frac{1}{2} \Big\{{\mathbb E}_{\mathbb P}[x_n(\delta)]-{\mathbb E}_{\mathbb P}[x_n(0)]\Big\}' \nabla^2_x g({\mathbb E}_{\mathbb P}[x_n(0)])\Big\{{\mathbb E}_{\mathbb P}[x_n(\delta)]-{\mathbb E}_{\mathbb P}[x_n(0)]\Big\}\\
& & + o\big(\|x_n(\delta)]-{\mathbb E}_{\mathbb P}[x_n(0)]\|^2\big) \\[5pt]
& = & g({\mathbb E}_{\mathbb P}[x_n(0)]) + \Big(\frac{1}{n}\,\nabla^2_x g(x^\star(0)) \overline{V}(0) + o(n^{-1})\Big)'\Big\{\delta \pi + O(\delta^2) + \frac{1}{n}O(\delta)\Big\} \\
& & + \frac{1}{2} \Big(\delta \pi + O(\delta^2) + \frac{1}{n} O(\delta)\Big)'\Big\{\nabla^2_x g(x^\star(0)) + O(n^{-1})\Big\}\Big(\delta \pi + O(\delta^2) + \frac{1}{n} O(\delta)\Big)
\end{eqnarray*}
where the second equality follows from \eqref{eq:expg_app} and \eqref{eq:exp-Ex-app}. It now follows that the expected reward at the mean of the robust solution
\begin{eqnarray}
\lefteqn{g({\mathbb E}_{\mathbb P}[x_n(\delta)])} \label{app:DRO-exp1} \\
& = & g({\mathbb E}_{\mathbb P}[x_n(0)]) + \underbrace{\frac{\delta}{n} \Big\{ \pi'\nabla^2_x g(x^\star(0))\overline{V}(0)\Big\} +  \frac{\delta^2}{2}\pi'\big[\nabla^2_x g(x^\star(0))\big]\pi  }_{\tiny \mbox{Impact of robustness}}\nonumber \\
& &  + o(1/n^2)+ o(\delta^2) + \frac{1}{n} o(\delta^2) + \frac{1}{n^2}O(\delta) \nonumber
\end{eqnarray}

In the case of the loss from Jensen's inequality, recall that the variance of the solution
\begin{eqnarray*}
\lefteqn{{\mathbb E}_{\mathbb P}\Big[\big(x_n(\delta))-{\mathbb E}_{\mathbb P}[x_n(\delta)]\big)\big(x_n(\delta))-{\mathbb E}_{\mathbb P}[x_n(\delta)]\big)'\Big]} \\
&   = & \frac{1}{n}\xi(\delta) + O(n^{-3/2}) \\ & = & \frac{1}{n}\Big\{\xi(0) + \delta \xi_\delta(0)\Big\} + \frac{1}{n}o(\delta) + O(n^{-3/2}).
\end{eqnarray*}
It now  follows from a Taylor series expansion that
\begin{eqnarray}
\lefteqn{{\mathbb E}_{\mathbb P} \Big[g(x_n(\delta))- g({\mathbb E}_{\mathbb P}[x_n(\delta)])\Big]} \nonumber \\
& = &  {\mathbb E}_{\mathbb P} \Big[\nabla_x g({\mathbb E}_{\mathbb P}[x_n(\delta)]) \Big(x_n(\delta)-{\mathbb E}_{\mathbb P}[x_n(\delta)]\Big)\Big] \nonumber \\
 & & +\frac{1}{2}{\mathbb E}_{\mathbb P} \Big[\big(x_n(\delta))-{\mathbb E}_{\mathbb P}[x_n(\delta)]\big)' \nabla^2_x g({\mathbb E}_{\mathbb P}[x_n(\delta)]) \big(x_n(\delta))-{\mathbb E}_{\mathbb P}[x_n(\delta)]\big)\Big] + O(n^{-3/2}) \nonumber \\
 & = &  
 \frac{1}{2}{\rm tr} \Big\{ \nabla^2_x g\Big(x^\star(\delta)+\frac{1}{n}\overline{V}(\delta)  + o(n^{-1}) \Big)  \frac{\xi(\delta)}{n}\Big\} + O(n^{-3/2}) \nonumber \\
 & = & \frac{1}{2n}\mbox{tr}\Big\{\nabla^2_x g(x^\star(\delta))\xi(\delta)\Big\}+ O(n^{-3/2})  \nonumber \\
 & = & \frac{1}{2n}\mbox{tr}\Big\{\nabla^2_x g(x^\star(0))\xi(0)\Big\}  + \frac{\delta}{2n}\frac{\mbox{d}}{\mbox{d}\delta}\Big\{\mbox{tr}\Big(\nabla^2_x g(x^\star(\delta))\xi(\delta)\Big)\Big\}\Big|_{\delta=0}+ \frac{1}{n}O(\delta^2) + O(n^{-3/2}) \nonumber \\
& = & {\mathbb E}_{\mathbb P} \Big[g(x_n(0))- g({\mathbb E}_{\mathbb P}[x_n(0)])\Big]+ \underbrace{\Big(\frac{\delta}{2n} \Big)\frac{\mbox{d}}{\mbox{d}\delta}\Big\{\mbox{tr}\Big(\nabla^2_x g(x^\star(\delta))\xi(\delta)\Big)\Big\}\Big|_{\delta=0}}_{\tiny \mbox{Impact of robustness}}
\label{app:DRO-exp2} \\
& & + \frac{1}{n}O(\delta^2) + O(n^{-3/2}) \nonumber
\end{eqnarray}
where the second equality follows from \eqref{eq:Ex-delta-app}, the fourth is a Taylor series expansion of the function $\mbox{tr}\Big\{\nabla^2_x g(x^\star(\delta))\xi(\delta)\Big\}$ around $\delta=0$, and
the last equality follows from \eqref{eq:decompTS2}. Proposition \ref{prop:out-of-sample reward} follows after substituting \eqref{app:DRO-exp1} and \eqref{app:DRO-exp2} into the decomposition \eqref{eq:DRO-decomp}.

\end{document}